\documentclass[10pt,reqno]{amsart}

\usepackage[utf8]{inputenc}
\usepackage{amsmath}
\usepackage{amssymb}
\usepackage{hyperref}
\usepackage{stackrel}
\usepackage{xcolor}
\usepackage{graphicx}
\usepackage{subfigure}
\usepackage{comment}
\usepackage{enumerate}
\usepackage{mathtools}

\newtheorem{theorem}{Theorem}[section]
\newtheorem{lemma}[theorem]{Lemma}

\newtheorem{definition}[theorem]{Definition}

\newtheorem{proposition}[theorem]{Proposition}
\newtheorem{remark}[theorem]{Remark}

\newtheorem{setting}[theorem]{Setting}

\DeclareMathOperator{\spann}{span}
\DeclareMathOperator{\conv}{conv}

\DeclareMathOperator{\ri}{ri}

\newcommand{\faceum}{A^\ast}
\newcommand{\facedois}{B^\ast}

\makeatletter
\newcommand{\tpitchfork}{%
  \vbox{
    \baselineskip\z@skip
    \lineskip-.52ex
    \lineskiplimit\maxdimen
    \m@th
    \ialign{##\crcr\hidewidth\smash{$-$}\hidewidth\crcr$\pitchfork$\crcr}
  }%
}
\makeatother

\makeatletter
\@namedef{subjclassname@2020}{\textup{2020} Mathematics Subject Classification}
\makeatother

\title{Smooth extremals on $C^0$-Finsler manifolds which are not locally minimizing}

\author{Ryuichi Fukuoka, J\'essica Buzatto Prudencio and Hugo Murilo Rodrigues}


\begin{document}

\maketitle

\begin{abstract}
Let $G$ be a three-dimensional Lie group, $\mathfrak{g}$ its Lie algebra and $\mathfrak{g}^\ast$ be the dual space of $\mathfrak{g}$. We consider left-invariant polyhedral Finsler structures on $G$ and study the local minimizing behavior of extremals using the Pontryagin Maximum Principle on the family of unit left-invariant vector fields of $G$.
In this setting, we present examples
on \(SE(2)\) and on the Heisenberg group \(\mathbb H\) for which a constant
Pontryagin extremal \((u_0,\mathfrak a_0) \in \mathfrak{g} \times \mathfrak{g}^{\ast}\) gives rise to a one-parameter
subgroup \(\gamma_{\mathfrak a_0}(t)=\exp(tu_0)\) that is not locally
minimizing. To prove this, we construct a one parameter family of competing extremals by concatenating left cosets of
one-parameter subgroups associated to the faces of the dual unit sphere  containing $\mathfrak{a}_0$.
In both Lie groups, these concatenated curves reach the same endpoints as
\(\gamma_{\mathfrak a_0}\), arbitrarily close to the identity, with strictly
smaller length. Thus, the examples show explicitly that an extremal,
even when associated with a constant vertical component and a unique constant maximizing
control, doesn't need to be locally minimizing.
This example is correlated to an asymptotic type curvature at $\mathfrak{a}_{0}$ defined on Finsler polyhedral structures on Lie groups.
\end{abstract}

\section{Introduction}
\label{section_introduction}

Let \(M\) be a smooth manifold, with tangent bundle \(TM\) and cotangent bundle
\(T^*M\). In this work, a \(C^0\)-Finsler structure on \(M\) is a continuous
function
\(
F:TM\longrightarrow [0,\infty)
\)
such that, for every \(x\in M\), the restriction
\(F_x:=F(x,\cdot):T_xM\to[0,\infty)\) is an asymmetric norm. We denote by
\(F_*:T^*M\to[0,\infty)\) its dual asymmetric norm, defined by
\(
F_*(x,\xi)=\max_{v\in S_{F_x}}\xi(v).
\)
The terminology for Finsler structures is not completely uniform in the
literature. Throughout this paper, the term \emph{Finsler structure}
will be reserved for the smooth case, whereas \(C^0\)-Finsler structure will
refer to the continuous asymmetric setting described above.

Geodesics and curvature are closely related in Riemannian and Finsler geometry.
One of the classical manifestations of this relation is provided by Jacobi
fields and conjugate points: curvature controls the infinitesimal spreading of
geodesics, whereas conjugate points provide an obstruction to their minimizing
behavior. In the nonsmooth setting of \(C^0\)-Finsler geometry, however, many
of the differential tools underlying this theory are no longer directly
available. Even the behavior of geodesics may differ substantially from the
smooth case: geodesics doesn't need to be smooth, and a prescribed initial position
and velocity may determine several geodesics or none at all (see \cite{Fukuoka-large-family}).

For a class of \(C^0\)-Finsler structures satisfying suitable horizontal
regularity assumptions, Pontryagin's Maximum Principle (PMP) provides an alternative
framework for studying minimizing curves. These structures were introduced in
\cite{Rodrigues} as \emph{Pontryagin-type} \(C^0\)-Finsler structures.
Let \(\{X_u\}_{u\in S^{n-1}}\) be the family of unit vector fields defining the
control system on an open subset $V$ of $M$ and, for \((x,\xi)\in T^*V\setminus\{0\}\), let
\[
\mathcal C(x,\xi)
=
\left\{
u\in S^{n-1}:
\xi(X_u(x))=\max_{v\in S^{n-1}}\xi(X_v(x))
\right\}.
\]
If \(H_u(x,\xi)=\xi(X_u(x))\) and \(\vec H_u\) denotes the Hamiltonian and the Hamiltonian vector
field on $T^\ast V\backslash 0$, the \emph{extended geodesic field} is the set-valued map
\[
\mathcal{E}(x,\xi)
=
\left\{
\vec H_u(x,\xi):u\in\mathcal C(x,\xi)
\right\}.
\]
The integral curves of the differential inclusion $(\dot{x}(t),\dot{\xi}(t)) \in \mathcal{E}(x(t), \xi(t))$ are Pontryagin extremals, and every length-minimizing curve on $V$
parametrized by arc length is  the projection $x(t)$ of a Pontryagin extremal.
Thus, PMP provides necessary conditions for minimality, but
they do not, in general, determine which extremals are locally minimizing.

Left-invariant \(C^0\)-Finsler structures on Lie groups provide a particularly
suitable framework in which to investigate this question. If \(G\) is a Lie
group, $\mathfrak{g}$ is its Lie algebra, and $\mathfrak{g}^{\ast}$ is de dual space of $\mathfrak{g}$, the Pontryagin extremal equation can be
represented on \(\mathfrak{g} \times \mathfrak g^*\) by the Euler--Arnold type equation \(\dot{\mathfrak a}(t)
= -\operatorname{ad}^{*}_{u(t)}\mathfrak a(t)\),
where \(u(t)\) is a maximizing control and \(\mathfrak a(t)\in\mathfrak g^*\)
is the vertical part of the extremal {see \cite{Prudencio-Fukuoka-1}}. For polyhedral Finsler structures, both
the unit sphere \(S_F\subset\mathfrak g\) and its dual sphere
\(S_{F_*}\subset\mathfrak g^*\) have a polyhedral structure that makes the
dynamics of the corresponding differential inclusion particularly explicit.

In \cite{Prudencio-Fukuoka-1}, this setting was studied for connected
three-dimensional Lie groups endowed with left invariant polyhedral Finsler
structures. The authors introduced an asymptotic expansion of flag curvature
by approximating a face of the polyhedral unit sphere by the unit spheres of a
family \(\{g_k\}_{k\in\mathbb N}\) of left-invariant Riemannian metrics. The
corresponding sectional curvatures are polynomial functions of \(k\), and the
dominant coefficient of this expansion defines the asymptotic curvature used
in the three-dimensional setting. This curvature was shown to be related to
the uniqueness of the maximizing control and, consequently, to the behavior of
extremals $u(t)$ corresponding to a given $\mathfrak{a}(t)$.

There is another reason why such an asymptotic curvature may be relevant.
First, recall the classical Riemannian picture. Let \(\gamma\) be a unit-speed
geodesic and suppose that, along \(\gamma\), the sectional curvature
\(K(\sigma)\) of every 2-plane \(\sigma\) containing \(\dot\gamma\) is
bounded from below by a positive constant \(k\). The comparison theorem for
Jacobi fields implies the occurrence of a conjugate point at a distance no
greater than \(\pi/\sqrt{k}\); moreover, a geodesic cannot remain minimizing
beyond its first conjugate point (see \cite{KobayashiNomizu2} and \cite{doCarmo2}). Consequently,
in the smooth case, larger positive curvature is associated with shorter
distance on which a geodesic can remain minimizing.

This classical picture suggests a possible interpretation of the asymptotic
curvature introduced in \cite{Prudencio-Fukuoka-1}. Indeed, this curvature
arises from sectional curvatures of Riemannian approximations whose unit
spheres converge locally to a face of the polyhedral unit sphere. Whenever the
dominant term of the corresponding curvature expansion is nonzero, the
sectional curvatures of these approximating metrics become unbounded as
\(k\to\infty\). From the Riemannian comparison picture, the associated
distance \(\pi/\sqrt{K}\) from which the extremal ceases to be minimizing, tends to zero. This suggests the following
heuristic principle: unbounded asymptotic curvature may be associated with an
obstruction to the local minimality in the
polyhedral geometry. We emphasize that, although the classical comparison
theorems cannot be applied directly to the nonsmooth limiting structure, they
provide the geometric intuition motivating the examples considered here.

The purpose of the present work is to provide concrete evidence for this
principle. We consider the three-dimensional Lie groups \(SE(2)\) (the special Euclidean group of the plane) and
\(\mathbb H\) (the Heisenberg
group), endowed with suitable left invariant polyhedral Finsler
structures. In each case, we construct a vertex
\(\mathfrak a_0\in S_{F_*}\) and a unique maximizing control \(u_0\), belonging
to the relative interior of the corresponding face in \(S_F\), such that
\((u_0,\mathfrak a_0)\) is a constant Pontryagin extremal. The corresponding
extremal on the group is the one-parameter subgroup
\(
\gamma_{\mathfrak a_0}(t)=\exp(tu_0).
\)
Although \(\gamma_{\mathfrak a_0}\) is an extremal, we show that it is not
locally minimizing.

Our proof is constructive. The Pontryagin extremals in a neighborhood of \(\mathfrak a_0\) are determined by periodic curves that resemble a center at \(\mathfrak a_0\). 
They are concatenations of smooth curves defined on the faces $\faceum$, $\facedois$ and $C^{\ast}$ of $S_{F_{\ast}}$ that contains \(\mathfrak a_0\). Their representation on the group are obtained by concatenating left cosets of one-parameter subgroups corresponding to the vertices $A, B, C \in S_F$, which are dual to the faces $A^\ast$, $B^\ast$ and $C^\ast$  respectively. We choose 
the traveling times in these concatenations so that the resulting curves meet 
$\gamma_{\mathfrak{a}_0}$ again arbitrarily close to the identity. A direct comparison 
then shows that the concatenated extremal has strictly smaller length than the 
corresponding segment of $\gamma_{\mathfrak{a}_{0}}$. Consequently, the constant extremal 
fails to be locally minimizing.

This work is organized as follows. In Section~\ref{prelimanaries}, we present the necessary preliminaries on asymmetric norms and dual asymmetric norms,
Pontryagin type \(C^0\)-Finsler structures, extended geodesic fields, and left-invariant structures on Lie groups. In Section~\ref{Milnors bases}, we present
Milnor bases for three-dimensional Lie algebras. 
Section~\ref{onthefaces} is devoted
to the behavior of the extended geodesic field on the relative interiors of
the faces of the dual unit sphere. In Section~\ref{transition between faces},
we study in the case when $\mathfrak{a}(t)$ are transversal
to the common edge. Finally, in
Section~\ref{in E(2) and H}, we construct extremals that are one-parameter
subgroups and are not locally minimizing. The cases of \(SE(2)\) and the
Heisenberg group are treated separately.

\section{Preliminaries}
\label{prelimanaries}

Convex analysis on a finite-dimensional real vector space \(V\) is often
developed after choosing an inner product \(\langle \cdot,\cdot\rangle\);
in particular, many presentations restrict the discussion to \(\mathbb{R}^n\)
endowed with its canonical inner product; see, for instance, \cite{Grunbaum} and \cite{RockafellarTyrrell}.
Once such an inner product is fixed, one obtains the usual identification
between \(V\) and its dual space \(V^*\), given by $v \in V \longmapsto \langle v,\cdot\rangle \in V^* .$
In the present work, however, the main objects are asymmetric norms $F:V\to \mathbb{R}$
and their dual norms $F_*:V^*\to \mathbb{R}.$
For this reason, it is more natural to formulate the discussion directly on
the dual space, using linear functionals \(\alpha\in V^*\), rather than
representing them through an inner product as \(\alpha=\langle v,\cdot\rangle\).

From this point to the end of the subsection, we shall use, with the necessary
adaptations to our setting, the theory  developed in \cite{Grunbaum}. 
We shall also use the terminology \emph{polyhedral set} instead of
\emph{polytope}.

\begin{definition}
Let \(V\) be a vector space. A map \(F:V\to [0,\infty)\) is called an
\emph{asymmetric norm} if it satisfies the following conditions:
\begin{enumerate}
    \item \(F(v)\geq 0\) for every \(v\in V\), and \(F(v)=0\) occurs only when
    \(v=0\);
    
    \item \(F(\lambda v)=\lambda F(v)\), for every \(v\in V\) and every
    \(\lambda>0\);
    
    \item \(F(v+w)\leq F(v)+F(w)\), for all \(v,w\in V\).
\end{enumerate}
\end{definition}

Asymmetric norms are generalizations of Minkowski norms (see
\cite{BaoChernShen}). Let \(V\) be endowed with an asymmetric norm \(F\). Given
\(y\in V\) and \(R>0\), we define the closed ball with center \(y\) and radius
\(R\) by
\[
    \overline{B}_{F}[y,R]
    :=
    \{z\in V \,;\, F(z-y)\leq R\}.
\]
Similarly, \(B_{F}(y,R)\) and \(S_{F}[y,R]\) will denote, respectively, the
open ball and the sphere with center \(y\) and radius \(R\). Since \(F\) is not
assumed to be symmetric, one should note that, in general,
\[
    F(z-y)\neq F(y-z).
\]
When the center is the origin and the radius is equal to one, we shall simply
write  $S_{F}$, $B_{F}$ and $\overline{B}_{F}$ for the unit sphere, the open unit ball and the closed unit ball, respectively.

\begin{definition}
Let \(F\) be an asymmetric norm on \(V\). The \emph{dual asymmetric norm}
associated with \(F\) is the function
\[
    F_*:V^*\to \mathbb{R}
\]
defined, for each \(\alpha\in V^*\), by
\[
    F_*(\alpha)
    :=
    \max_{u\in \overline{B}_{F}} \alpha(u)=\max_{u\in S_F} \alpha(u).
\]

\end{definition}

Analogously, given \(\alpha_0\in V^*\) and \(R>0\), we denote  the closed ball, the open ball
and the sphere associated with the dual asymmetric norm \(F_*\), respectively,
by  $\overline{B}_{F_*}[\alpha_0,R]$, $B_{F_*}(\alpha_0,R)$ and $S_{F_*}[\alpha_0,R]$ respectively.

When the center is the origin and the radius is equal to one, we shall use the
simplified notation $\overline{B}_{F_*}$, $B_{F_*}$ and $S_{F_*}$ for the closed unit ball, the open unit ball and the unit sphere of \(F_*\), respectively.

A subset of \(V\) will be called \emph{polyhedral} when it can be written as an intersection of closed half-spaces of \(V\).

\begin{definition}
Let \(F:V\to \mathbb{R}\) be an asymmetric norm. We say that \(F\) is a
\emph{polyhedral asymmetric norm} if its closed unit ball $\overline{B}_{F}$ is a polyhedral subset of \(V\).
\end{definition}

\subsection{Extended geodesic fields}

In this subsection, we recall the notion of a Pontryagin-type
\(C^0\)-Finsler   structure on a smooth \(n\)-dimensional manifold \(M\). This
class of structures, introduced in \cite{Rodrigues} (see also \cite{Fukuoka-Rodrigues}), is designed to satisfy the
minimal regularity assumptions needed in the framework of Pontryagin's
maximum principle. In this setting, the control set is taken to be the
Euclidean unit sphere $\mathcal{C}=S^{n-1}$, and the control system is given by a family of $C^1$ unit vector fields defined on an open subset $U\subset M$.
The construction leads to a structure that includes the geodesic spray of a Finsler manifold, and is called
\emph{extended geodesic field}: it is a multivalued vector field defined on the slit cotangent bundle  $T^*U\setminus\{0\}$.

\begin{definition}
  A \emph{\(C^0\)-Finsler structure} on \(M\) is
a continuous function $F:TM\to \mathbb{R}$ 
such that, for every \(x\in M\), the restriction $F_x:=F(x,\cdot):T_xM\to \mathbb{R}$
is an asymmetric norm on the vector space \(T_xM\).
\end{definition}

For each \(x\in M\), we shall denote by $S_{F_x}, B_{F_x}$ and $\overline{B}_{F_x}$
the unit sphere, the open unit ball and the closed unit ball associated with
the asymmetric norm \(F_x\), all centered at the origin of \(T_xM\),
respectively.

\begin{definition}\label{Pontryagin type C0-Finsler manifold}
A \(C^0\)-Finsler manifold is of Pontryagin type at \(p \in M\) if there exist
a neighborhood \(U\) of \(p\), a coordinate system
\(\phi = (x^1, \ldots, x^n): U \to \mathbb{R}^n\)
(with the respective natural coordinate system
\(\phi_{TU} := d\phi : (x^1, \ldots, x^n, y^1, \ldots, y^n): TU \to
\mathbb{R}^{2n}\) of the tangent bundle), and a family of \(C^1\) unit vector
fields
\[
\left\{
x \mapsto X_u(x)
=
\left(
y^1(x^1,\ldots,x^n,u), \ldots,
y^n(x^1,\ldots,x^n,u)
\right);
\ u \in S^{n-1}
\right\}
\]
on \(U\) parametrized by \(u \in S^{n-1}\) such that
\begin{enumerate}
    \item \(u \mapsto X_u(x)\) is a homeomorphism from
    \(S^{n-1} \subset \mathbb{R}^n\) onto \(S_{F_x}\) for every
    \(x \in U\);

    \item \((x,u) \mapsto
    \left(
    y^1(x^1,\ldots,x^n,u), \ldots,
    y^n(x^1,\ldots,x^n,u)
    \right)\) is continuous;

    \item \((x,u) \mapsto
    \left(
    \dfrac{\partial y^1}{\partial x^i}(x^1,\ldots,x^n,u),
    \ldots,
    \dfrac{\partial y^n}{\partial x^i}(x^1,\ldots,x^n,u)
    \right)\) is continuous for every
    \(i = 1,\ldots,n\).
\end{enumerate}
We say that \(F\) is of Pontryagin type on \(M\) if \((M,F)\) is of
Pontryagin type at every \(p \in M\). In this case \((M,F)\) is a Pontryagin-type \(C^0\)-Finsler manifold.
\end{definition}

The notion introduced in Definition \ref{Pontryagin type C0-Finsler manifold}
does not depend on the particular system of coordinates chosen. More generally,
it is also independent of the coordinate system on \(TU\) induced by a local
trivialization (see \cite{Rodrigues}, \cite{Fukuoka-Rodrigues}).

Let us write
\[
    X_u(x)
    =
    f^i(x,u)\frac{\partial}{\partial x^i}
    =
    \bigl(f^1(x,u),\ldots,f^n(x,u)\bigr),
\]
where, from now on, we use the Einstein summation convention. For each
\(u\in S^{n-1}\), consider the Hamiltonian $ H_u:T^*U\to\mathbb{R}$
defined by $ H_u(x,\xi)=\xi\bigl(X_u(x)\bigr).$
The Hamiltonian vector field associated with \(H_u\) will be denoted by
\(\vec{H}_u\), and is characterized by $
    dH_u(x,\xi)
    =
    \omega\bigl(\,\cdot\,,\vec{H}_u(x,\xi)\bigr),$
    where $\omega$ is the symplectic 2-form.
The \emph{extended geodesic field} is then defined as the set-valued map
\(\mathcal{E}\) which assigns to each \((x,\xi)\in T^*U\setminus\{0\}\) the set
\begin{equation}
    \mathcal{E}(x,\xi)
    =
    \left\{
        \vec{H}_u(x,\xi)\,;\,
        u\in \mathcal{C}(x,\xi)
    \right\}.
\end{equation}
In the natural coordinates $(x,\xi_i dx^i)\longmapsto
    (x^1,\ldots,x^n,\xi_1,\ldots,\xi_n)$
on \(T^*U\), this set can be written as
\begin{equation}
    \mathcal{E}(x,\xi)
    =
    \left\{
        f^i(x,u)\frac{\partial}{\partial x^i}
        -
        \xi_j
        \frac{\partial f^j}{\partial x^i}(x,u)
        \frac{\partial}{\partial \xi_i}
        \,;\,
        u\in \mathcal{C}(x,\xi)
    \right\},
\end{equation}
where $\mathcal{C}(x,\xi)
    =
    \left\{
        u\in S^{n-1}\,;\,
        \xi\bigl(X_u(x)\bigr)
        =
        \max_{v\in S^{n-1}}
        \xi\bigl(X_v(x)\bigr)
    \right\}.$
An absolutely continuous curve $(x(t),\xi(t)):[a,b]\to T^*U$
is called a \emph{Pontryagin extremal} of \((U,F)\) if
\begin{equation}
    \bigl(\dot{x}(t),\dot{\xi}(t)\bigr)
    \in
    \mathcal{E}\bigl(x(t),\xi(t)\bigr)
\end{equation}
for almost every \(t\in [a,b]\). In this case, the projected curve \(x(t)\) is
called an \emph{extremal} of \((U,F)\).

By Pontryagin's Maximum Principle, if \(x(t)\) is a length-minimizing curve in
\(U\), then there exists a Pontryagin extremal \((x(t),\xi(t))\) of \((U,F)\)
such that
\begin{equation}
    \mathcal{M}\bigl(x(t),\xi(t)\bigr)
    =
    \max_{u\in S^{n-1}}
    \xi(t)\bigl(X_u(x(t))\bigr)
    =
    \max_{v\in S_{F_{x(t)}}}
    \xi(t)(v)
    =
    (F_*)_{x(t)}\bigl(\xi(t)\bigr)   
\end{equation}
is constant along the curve.

For a smooth Finsler manifold $(M,F)$, the fundamental tensor gives rise, fibre by fibre, to the Legendre transform  $\mathcal{L}$ between $TM\setminus 0$ and $T^{*}M\setminus 0$, which is a smooth map. 
 In \cite{Rodrigues}, \cite{Fukuoka-Rodrigues}, it was shown that  $\mathcal{E}=d\mathcal{L}(\mathbf{G})$, where $\mathbf{G}$ is the geodesic spray of $(M,F)$.

\subsection{Left-invariant \texorpdfstring{\(C^0\)}{C⁰}-Finsler structures on Lie groups}

\begin{definition}
Let \(G\) be a Lie group and let \(e\in G\) denote its identity element. For
each \(g\in G\), let \(L_g:G\to G\) be the left translation by \(g\). A
\(C^0\)-Finsler structure \(F\) on \(G\) is said to be \emph{left invariant}
if
\[
    F(x,y)=F\bigl(gx,d(L_g)_x(y)\bigr),
\]
for every \(g,x\in G\) and every \(y\in T_xG\).
\end{definition}

\begin{remark}
There is a bijective correspondence between left invariant \(C^0\)-Finsler
structures on \(G\) and asymmetric norms on its Lie algebra
\(\mathfrak{g}=T_eG\). Indeed, as in the Riemannian case, the structure is
completely determined by its restriction to \(T_eG\), and the values on the
other tangent spaces are obtained by left translations.
\end{remark}

 On a Lie group, mathematical objects on $G$ can be represented on $\mathfrak{g}$ and $\mathfrak{g}^\ast$. In Theorem \ref{infinitesimal geodesic}, we describe the extended geodesic field of \((G,F)\) on
\(S_{F_e}\times\mathfrak{g}^*\). In this formulation, the control set
\(S^{n-1}\) is replaced by the unit sphere \(S_{F_e}\subset\mathfrak{g}\),
which is more convenient for the problem considered here. 
 The control system is given by the left-invariant unit vector fields $X_u$ such that $X_u(e)=u$. This result  can be found in   \cite{Prudencio-Fukuoka-1}.

\begin{theorem}\label{infinitesimal geodesic}
Let \((G,F)\) be a Lie group endowed with a left invariant \(C^0\)-Finsler
structure. Let \(t\in I\mapsto (x(t),\xi(t))\) be a Pontryagin extremal of
\((G,F)\) associated with a control \(u(t)\). Define
\(\mathfrak{a}(t)=(dL_{x(t)})^*_{e}(\xi(t))\in\mathfrak{g}^*\). Then \(\mathfrak{a}(t)\) is
absolutely continuous and \(u(t)\in\mathcal{C}(\mathfrak{a}(t))\) is a measurable control
satisfying
\begin{equation}\label{infinitesimal cotangent form}
    \dot{\mathfrak{a}}(t)
    =
    -\operatorname{ad}^{*}(u(t))(\mathfrak{a}(t))
    =
    \mathfrak{a}(t)\bigl([u(t),\cdot]\bigr).
\end{equation}

Conversely, suppose that \(\mathfrak{a}(t)\in\mathfrak{g}^*\setminus\{0\}\) is absolutely
continuous and that \(u(t)\in\mathcal{C}(\mathfrak{a}(t))\) is a measurable control
satisfying \eqref{infinitesimal cotangent form}. Then, for every \(t_0\in I\) and every \(x_0\in G\), there
exists a unique Pontryagin extremal \((x(t),\xi(t))\) associated with the
control \(u(t)\) such that
\(\mathfrak{a}(t)=(dL_{x(t)})^*_{e}(\xi(t))\) and \(x(t_0)=x_0\).
\end{theorem}

In this work, we denote by \(\mathcal{E}_{\mathfrak{g}^*}\) the set-valued map
which assigns to each \(\mathfrak a\in\mathfrak{g}^*\setminus\{0\}\) the set
\[
    \mathcal{E}_{\mathfrak{g}^*}(\mathfrak a)
    =
    \left\{
        -\operatorname{ad}^{*}(u)(\mathfrak a)
        \,;\,
        u\in\mathcal{C}(\mathfrak a)
    \right\}.
\]
Therefore, every solution \((u(t),\mathfrak{a}(t))\) of \eqref{infinitesimal cotangent form} determines a curve
\(\mathfrak{a}(t)\) satisfying the differential inclusion
\begin{equation}\label{inclusão}
    \dot{\mathfrak{a}}(t)\in\mathcal{E}_{\mathfrak{g}^*}\bigl(\mathfrak{a}(t)\bigr).
\end{equation}

\section{Milnor's bases for three-dimensional Lie algebras}
\label{Milnors bases}

In this section, we recall the notion of Milnor bases for
three-dimensional Lie algebras \(\mathfrak{g}\)  (see \cite{Milnor}). For each Lie
algebra considered in this work, one such basis will be fixed.

\begin{theorem}
Let \(G\) be a connected and oriented three-dimensional Lie group, and let
\(\mathfrak{g}\) be its Lie algebra endowed with a left invariant Riemannian
metric \(\langle\cdot,\cdot\rangle\). Denote by
\((u,v)\mapsto u\times v\) the cross product on \(\mathfrak{g}\) determined by
\(\langle\cdot,\cdot\rangle\) and by the chosen orientation. Then there exists
a unique linear map \(A:\mathfrak{g}\to\mathfrak{g}\) such that
\[
    [u,v]=A(u\times v)
\]
for all \(u,v\in\mathfrak{g}\). Moreover, the Lie group \(G\) is unimodular if
and only if \(A\) is self-adjoint.
\end{theorem}

\begin{proof}
See \cite{Milnor}.
\end{proof}

\begin{setting}\label{Milnor setting}
Throughout this work, for each Lie algebra \(\mathfrak{g}\), we fix an
auxiliary basis
\(\mathcal{B}=\{e_1,e_2,e_3\}\). This basis determine both an orientation and
an inner product \(\langle\cdot,\cdot\rangle\) on \(\mathfrak{g}\), chosen so
that \(\mathcal{B}\) is a positive orthonormal basis.

These choices also induce an orientation and an inner product on the dual
space \(\mathfrak{g}^*\), in such a way that the dual basis
\(\mathcal{B}^*=\{e^1,e^2,e^3\}\) is positive and orthonormal. We shall use the
musical isomorphisms
\(\flat:\mathfrak{g}\to\mathfrak{g}^*\) and
\(\sharp:\mathfrak{g}^*\to\mathfrak{g}\) induced by
\(\langle\cdot,\cdot\rangle\). Thus, elements of \(\mathfrak{g}\) may be
identified with elements of \(\mathfrak{g}^*\) having the same coordinates
with respect to the bases \(\mathcal{B}\) and \(\mathcal{B}^*\).

Finally, the bases \(\mathcal{B}\) and \(\mathcal{B}^*\) determine cross
products on \(\mathfrak{g}\) and on \(\mathfrak{g}^*\), respectively. These
will be denoted by \(\times_{\mathfrak{g}}\) and
\(\times_{\mathfrak{g}^*}\); when no confusion is possible, both will simply be
written as \(\times\).
\end{setting}

If \(\mathfrak{g}\) is a unimodular tridimensional Lie algebra, we can find a
basis \(\tilde{\mathcal{B}}=\{\tilde{e}_1,\tilde{e}_2,\tilde{e}_3\}\) such
that
\[
    [\tilde{e}_1,\tilde{e}_2]=\tilde{\lambda}_3\tilde{e}_3,\qquad
    [\tilde{e}_2,\tilde{e}_3]=\tilde{\lambda}_1\tilde{e}_1
    \quad\text{and}\quad
    [\tilde{e}_3,\tilde{e}_1]=\tilde{\lambda}_2\tilde{e}_2,
\]
what means that
\[
    [A]_{\tilde{\mathcal{B}}}
    =
    \begin{bmatrix}
        \tilde{\lambda}_1 & 0 & 0 \\
        0 & \tilde{\lambda}_2 & 0 \\
        0 & 0 & \tilde{\lambda}_3
    \end{bmatrix}
\]
if we consider \(\tilde{\mathcal{B}}\) as the positive orthonormal basis, see
\cite{Milnor}. We will modify this basis (and the corresponding inner product) to get \(\breve{\mathcal{B}}\).

Set \(\operatorname{sgn}(\tilde{\lambda})=1\) if \(\tilde{\lambda}>0\),
\(\operatorname{sgn}(\tilde{\lambda})=0\) if \(\tilde{\lambda}=0\) and
\(\operatorname{sgn}(\tilde{\lambda})=-1\) if \(\tilde{\lambda}<0\). Define
\[
    \breve{e}_i
    =
    \begin{cases}
        \tilde{e}_i,
        & \text{if } \tilde{\lambda}_j=\tilde{\lambda}_k=0; \\[0.3cm]
        \dfrac{\tilde{e}_i}{\sqrt{|\tilde{\lambda}_j|}},
        & \text{if } \tilde{\lambda}_j\neq 0
        \text{ and } \tilde{\lambda}_k=0; \\[0.5cm]
        \dfrac{\tilde{e}_i}{\sqrt{|\tilde{\lambda}_j\tilde{\lambda}_k|}},
        & \text{if } \tilde{\lambda}_j\neq 0
        \text{ and } \tilde{\lambda}_k\neq 0,
    \end{cases}
\]
where \(\{i,j,k\}=\{1,2,3\}\). Then
\(\breve{\mathcal{B}}=\{\breve{e}_1,\breve{e}_2,\breve{e}_3\}\) is a basis of
\(\mathfrak{g}\) such that
\[
    [\breve{e}_1,\breve{e}_2]
    =
    \operatorname{sgn}(\tilde{\lambda}_3)\breve{e}_3,\qquad
    [\breve{e}_2,\breve{e}_3]
    =
    \operatorname{sgn}(\tilde{\lambda}_1)\breve{e}_1
\]
and
\[
    [\breve{e}_3,\breve{e}_1]
    =
    \operatorname{sgn}(\tilde{\lambda}_2)\breve{e}_2.
\]
Set \(\lambda_i=\operatorname{sgn}(\tilde{\lambda}_i)\). Therefore,
\[
    [A]_{\breve{\mathcal{B}}}
    =
    \begin{bmatrix}
        \lambda_1 & 0 & 0 \\
        0 & \lambda_2 & 0 \\
        0 & 0 & \lambda_3
    \end{bmatrix}
\]
with respect to the positive orthonormal basis \(\breve{\mathcal{B}}\).
Finally order \(\breve{e}_1,\breve{e}_2\) and \(\breve{e}_3\) in such a way
that \(\lambda_i=0\) comes first, \(\lambda_i=-1\) comes right after and
\(\lambda_i=1\) comes last. This reordering of \(\breve{\mathcal{B}}\) give us
the ordered orthonormal basis
\(\mathcal{B}=\{e_1,e_2,e_3\}\) for the unimodular case. With this setup, the
unimodular Lie algebras are the following  (see \cite{Milnor}): 
\begin{enumerate}
    \item \(\mathfrak{a}\) (Abelian algebra):
    \((\lambda_1,\lambda_2,\lambda_3)=(0,0,0)\);

    \item \(\mathfrak{h}\) (Heisenberg algebra):
    \((\lambda_1,\lambda_2,\lambda_3)=(0,0,1)\);

    \item \(\mathfrak{e}(2)\):
    \((\lambda_1,\lambda_2,\lambda_3)=(0,1,1)\);

    \item \(\mathfrak{e}(1,1)\):
    \((\lambda_1,\lambda_2,\lambda_3)=(0,-1,1)\);

    \item \(\mathfrak{so}(1,2)\):
    \((\lambda_1,\lambda_2,\lambda_3)=(-1,1,1)\);

    \item \(\mathfrak{so}(3)\):
    \((\lambda_1,\lambda_2,\lambda_3)=(1,1,1)\).
\end{enumerate}

\begin{remark}
In this work, we are interested in studying extremals in $SE(2)$ and $\mathbb{H}$ considering $\mathcal{E}_{\mathfrak{g}^{\ast}}$ on $\mathfrak{h}$ and $\mathfrak{e}(2)$ .
\end{remark}

\begin{remark}
Unless otherwise stated, we are in Settings \ref{Milnor setting} throughout this work.
\end{remark}

\section{\texorpdfstring{$\mathcal{E}_{\mathfrak{g}^*}$}{E(g*)} on the faces of \texorpdfstring{$S_{F_*}$}{S(F*)}.}
\label{onthefaces}

From now on we will consider just three-dimensional Lie groups and in order to simplify the notation, we shall denote the unit sphere of \((\mathfrak{g},F_e)\) by \(S_F\), and the unit sphere of
\((\mathfrak{g}^*,F_*)\) by \(S_{F_*}\). This section is devoted to describing
\(\mathcal{E}_{\mathfrak{g}^{\ast}}\vert_{\ri H}\), where $\ri H$ is the relative interior of a face $H$ of $S_{F_{\ast}}$.

Let \(G\) be a Lie group endowed with a left invariant polyhedral Finsler
structure \(F\). Recall that the extended geodesic field
\(\mathcal{E}_{\mathfrak{g}^{\ast}}\) is given by
\[
    \mathcal{E}_{\mathfrak{g}^{\ast}}(\mathfrak{a})
    =
    \left\{
        \mathfrak{a}([u,\cdot])\in \mathfrak{g}^*
        \,;\,
        u\in \mathcal{C}(\mathfrak{a})
    \right\},
    \qquad \mathfrak{a}\in S_{F_*}.
\]

\begin{proposition}
Let \(H\) be a face of \(S_{F_*}\), and let \(u_H\in S_F\) be the unique
vertex that maximizes every \(\mathfrak a\in\ri H\). Then
\(\mathcal{E}_{\mathfrak{g}^*}\vert_{\ri H}\) is the smooth vector field
\[
    \mathfrak{a}\in {\ri} H
    \longmapsto
    X_{u_H}(\mathfrak{a}) 
    =
    \mathfrak{a}([u_H,\cdot]). 
\]
\end{proposition}
\begin{proof}
It suffices to note that, under these conditions, for each chosen $\mathfrak{a}$, there exists a unique $u$ that maximizes it.
\end{proof}

 Let \(\mathcal{B}\) be the basis of \(\mathfrak{g}\) fixed in
Settings \ref{Milnor setting}. We denote by \(A^T:\mathfrak{g}^*\to\mathfrak{g}^*\) the
transpose of \(A:\mathfrak{g}\to\mathfrak{g}\), and by
\(A^*:\mathfrak{g}\to\mathfrak{g}\) its adjoint with respect to the inner
product induced by \(\mathcal{B}\). From the definitions of the musical
isomorphisms, it is straightforward that
\[
    \bigl(A^T(\mathfrak{a})\bigr)^\sharp=A^*(\mathfrak{a}^\sharp).
\]

\begin{proposition}
Let \(\mathcal{B}\) and \(\mathcal{B}^*\) be the bases of
\(\mathfrak{g}\) and \(\mathfrak{g}^*\), respectively, introduced in
Settings \ref{Milnor setting}. Then
\begin{equation}\label{Xu by AT x}
    X_u(\mathfrak{a})=A^T(\mathfrak{a})\times u^\flat .
\end{equation}
\end{proposition}

\begin{proof}
Let \(v\in\mathfrak{g}\). Using the definition of \(X_u\), the relation
\([u,v]=A(u\times v)\), and the definition of the transpose map, we obtain
\[
    X_u(\mathfrak{a})(v)
    =
    \mathfrak{a}([u,v])
    =
    \mathfrak{a}(A(u\times v))
    =
    A^T(\mathfrak{a})(u\times v).
\]
Identifying \(\mathfrak{g}\) and \(\mathfrak{g}^*\) through the musical
isomorphisms, this last term can be written as
\[
    A^T(\mathfrak{a})(u\times v)
    =
    \left\langle A^*(\mathfrak{a}^\sharp),u\times v\right\rangle.
\]
By the properties of the mixed product,
\[
    \left\langle A^*(\mathfrak{a}^\sharp),u\times v\right\rangle
    =
    \left\langle A^*(\mathfrak{a}^\sharp)\times u,v\right\rangle.
\]
Since \(\bigl(A^T(\mathfrak{a})\bigr)^\sharp=A^*(\mathfrak{a}^\sharp)\), this is precisely $\bigl(A^T(\mathfrak{a})\times u^\flat\bigr)(v).$ Therefore,
\[
    X_u(\mathfrak{a})=A^T(\mathfrak{a})\times u^\flat,
\]
as claimed.
\end{proof}

We next write the matrix form of \(X_u\) with respect to the bases chosen in
Settings \ref{Milnor setting}. This computation follows directly either from the previous
identity or from the expression \(X_u=\mathfrak{a}([u,\cdot])\). Thus, if
\(\mathfrak{a}=(\mathfrak{a}_1,\mathfrak{a}_2,\mathfrak{a}_3)\) and
\(u=(u^1,u^2,u^3)\), then
\begin{equation}\label{fildsbymatrix}
    X_u(\mathfrak{a})
    =
    [X_u][\mathfrak{a}]
    =
    \begin{bmatrix}
        0 & \lambda_2 u^3
          & -\lambda_3 u^2 \\[0.15cm]
        0 & 0
          & \lambda_3 u^1 \\[0.15cm]
        0 & -\lambda_2 u^1
          & 0
    \end{bmatrix}
    \begin{bmatrix}
        \mathfrak{a}_1\\
        \mathfrak{a}_2\\
        \mathfrak{a}_3
    \end{bmatrix}.
\end{equation}
\section{Transition of \texorpdfstring{$\mathcal{E}_{\mathfrak{g}^*}$}{E(g*)} between faces.}
\label{transition between faces}
In the previous section, \eqref{fildsbymatrix} gave the coordinate expression of the vector field \(\mathcal{E}_{\mathfrak{g}^*}|_{\ri H}\) once the maximizing vertex \(u_H\) is fixed. In this section, we investigate how this vector field behaves when it reaches an edge and how the transition of $\mathfrak{a}(t)$ from one face to another  occur when its incidence on the edge is transversal.

Let \( H \) be a face of \( S_{F_*} \) and \( E \) be an edge contained in \( H \). Let \( \eta \) be a unit functional  parallel to \( E \), \( n \) be a unit functional parallel to \( H \) and orthogonal to \( E \). Then
\begin{equation}\label{uflat n eta}
\mathcal{B} = \left( \frac{u_H^\flat}{\|u_H^\flat\|}, n, \eta \right)
\end{equation}
is an orthonormal basis of \( \mathfrak{g}^* \). Set \( \sigma \left( u_H^\flat, n, \eta \right) = 1 \) if the basis is positive and \( \sigma \left( u_H^\flat, n, \eta \right) = -1 \) otherwise.

\begin{proposition}
Let \( E \) be an edge contained in \( H \) and \(\mathfrak{a}_0 \in E \). Parametrize \( \mathfrak{a} \) by \( t \mapsto  \mathfrak{a}_0 + t\eta \). Then
\begin{equation}\label{scalar Xu by AT}
\sigma \left( u^\flat, n, \eta \right) \frac{\langle X_u(\mathfrak{a}_0 + t\eta), n \rangle}{\|u^\flat\|} = \langle A^T(\mathfrak{a}_0), \eta \rangle + t \langle A^T(\eta), \eta \rangle.
\end{equation}
In particular
\begin{equation}\label{scalar Xu by AT t=0}
\sigma \left( u^\flat, n, \eta \right) \frac{\langle X_u(\mathfrak{a}_0), n \rangle}{\|u^\flat\|} = \langle A^T(\mathfrak{a}_0), \eta \rangle.
\end{equation}
\end{proposition}

\begin{proof}
We have
\[
\begin{aligned}
&\frac{\langle X_u(\mathfrak{a} + t\eta), n \rangle}{\|u^\flat\|} = \frac{\langle A^T(\mathfrak{a} + t\eta) \times u^\flat, n \rangle}{\|u^\flat\|} 
= \frac{\langle A^T(\mathfrak{a}) \times u^\flat, n \rangle}{\|u^\flat\|} + t \frac{\langle A^T(\eta) \times u^\flat, n \rangle}{\|u^\flat\|} \\
&= \langle A^T(\mathfrak{a}), \frac{u^\flat}{\|u^\flat\|} \times n \rangle + t \langle A^T(\eta), \frac{u^\flat}{\|u^\flat\|} \times n \rangle \\
&= \sigma \left( u^\flat, n, \eta \right) \left( \langle A^T(\mathfrak{a}), \eta \rangle + t \langle A^T(\eta), \eta \rangle \right),
\end{aligned}
\]
where the first equation is due to (\ref{Xu by AT x}), the third equation is a symmetry of the mixed product and the last equation is due to (\ref{uflat n eta}).
\end{proof}

\begin{setting}\label{setting to edge}
Let \( E \) be an edge of \( S_{F_*} \). Denote the faces that contain \( E \) by \( \hat{H} \) and \( \tilde{H} \) and set \( \hat{u} = u_{\hat{H}} \) and \( \tilde{u} = u_{\tilde{H}} \). Notice that they are the vertices of \( S_F \) that are maximized by \( \operatorname{ri} \hat{H} \) and \( \operatorname{ri} \tilde{H} \) respectively. Then \( \operatorname{ri} E \) is maximized by the edge $\operatorname{conv}\{\hat{u}, \tilde{u}\} = \{(1 - \theta)\hat{u} + \theta \tilde{u} : \theta \in [0, 1]\} \subset S_F.$
Define \( \hat{n} \in \mathfrak{g}^{\ast} \) as the \( \|\cdot\| \)-unit vector parallel to \( \hat{H} \), orthogonal to \( E \) and pointed outward from \( \hat{H} \). 
Define \( \tilde{n} \) with the same properties of \(\hat{n}\) but with respect to \(\tilde{H}\). Let \(\eta = (\eta_1, \eta_2, \eta_3)\) be a \(\|\cdot\|\)-unit functional parallel to \(E\) such that
$\hat{\mathcal{B}} = \{\hat{u}^\flat, \hat{n}, \eta\}$
is oriented positively. Observe that
$\tilde{\mathcal{B}} = \{\tilde{u}^\flat, \tilde{n}, \eta\}$
is oriented negatively.
\end{setting}

\begin{proposition}\label{flow corresponding}
Suppose that we are in Settings (\ref{setting to edge}). Then
\begin{equation}\label{flow corresponding eq}
\frac{\langle X_{\hat{u}}(\mathfrak{a}), \hat{n} \rangle}{\|\hat{u}^\flat\|} = -\frac{\langle X_{\tilde{u}}(\mathfrak{a}), \tilde{n} \rangle}{\|\tilde{u}^\flat\|}.
\end{equation}
 
\end{proposition}

\begin{proof}
Both sides are equal to \(\langle A^T(\mathfrak{a}), \eta \rangle\) due to (\ref{scalar Xu by AT t=0}).
\end{proof}

\begin{remark}
Proposition \ref{flow corresponding} states that \(X_{\tilde{u}}\) ``enters transversally'' into \(\tilde{H}\), iff \(X_{\hat{u}}\) ``leaves transversally'' \(\hat{H}\). Moreover, the flow corresponding to \(X_{\hat{u}}/\|\hat{u}^\flat\|\) that ``enters'' into \(\tilde{H}\) is exactly the flow corresponding to \(X_{\hat{u}}/\|\hat{u}^\flat\|\) that leaves \(\hat{H}\).  

Set  
\begin{equation}\label{T_E}
\mathcal{T}_E := \{\mathfrak{a} \in \operatorname{ri} E; X_{\tilde{u}}(\mathfrak{a}) \text{ is tangent to } \operatorname{ri} E\}.
\end{equation}

Notice that we also have  
\[
\mathcal{T}_E := \{\mathfrak{a} \in \operatorname{ri} E; X_{\hat{u}}(\mathfrak{a}) \text{ is tangent to } \operatorname{ri} E\}
\]
due to Proposition \ref{flow corresponding}.
\end{remark}

\begin{definition}
The set \(\mathcal{T}_E\) given by (\ref{T_E}) is the set of tangency points of \(\mathcal{E}_{\mathfrak{g}^*}\) in \(\operatorname{ri} E\).
\end{definition}

For our purposes, it will be sufficient to understand what happens to \(\mathcal{E}_{\mathfrak{g}^*}\) in a neighborhood of $ \mathfrak{a}_0\in \operatorname{ri} E\backslash \mathcal{T}_E$.
In this case, an integral curve of \(\mathcal{E}_{\mathfrak{g}^*}\) in a neighborhood of \(\mathfrak{a}_0\) is obtained by concatenating integral curves on \(\tilde H\) and \(\hat H\) (or vice versa).

\begin{lemma}\label{across the edge}
Suppose that we are in Setting \ref{setting to edge}. Let \(\mathfrak a(t)\) be the vertical part of a Pontryagin extremal  such that \( \mathfrak{a}(t_0) \in \operatorname{ri} E \setminus \mathcal{T}_E \). Without loss of generality, suppose that \( \langle X_{\tilde{u}}(\mathfrak{a}(t_0)), \tilde{n} \rangle > 0 \). Then there exist \( \varepsilon > 0 \) such that \( \mathfrak{a}((t_0 - \varepsilon, t_0)) \) lies in \( \operatorname{ri} \tilde{H} \) and \( \mathfrak{a}((t_0, t_0 + \varepsilon)) \) lies in \( \operatorname{ri} \hat{H} \).
\end{lemma}

\begin{proof}
Consider the signed distance function \( \sigma : \tilde{H} \cup \hat{H} \to \mathbb{R} \) defined by
\[
\sigma(\mathfrak{a}) = 
\begin{cases} 
-\operatorname{dist}_{\|\cdot\|}(\mathfrak{a}, r) & \text{if } \mathfrak{a} \in \tilde{H}; \\ 
\operatorname{dist}_{\Vert \cdot\Vert}(\mathfrak{a}, r) & \text{if } \mathfrak{a} \in \hat{H}, 
\end{cases}
\]
where \( r \) is the line containing \( E \). Its level sets are parallel to \( E \) and \( \sigma(\mathfrak{a}(t)) \) is nondecreasing function in a neighborhood \( J \) of \( t_0 \) due to \( \langle X_{\tilde{u}}(\mathfrak{a}(t_0)), \tilde{n} \rangle > 0 \) and  (\ref{flow corresponding eq}). 
Notice that if \( t_1 \) is another point in \( J \) such that \( \sigma(\mathfrak{a}(t_1)) = 0 \), then \( \sigma \circ \mathfrak{a}|_{[t_0, t_1]} = 0 \) (or else \( \sigma \circ \mathfrak{a}|_{[t_1, t_0]} = 0 \)), which would imply that \( \mathfrak{a}|_{[t_0, t_1]} \in E \).  We will prove that this can't happen.

 Recall that if $\mathfrak{a}_0 \in \ri E$ and $E^\sharp\subset S_{F}$ is the dual edge of $E$, then $\mathfrak{a}_0$ is maximized by $\mathcal{C}(\mathfrak{a}_0) = E^\sharp = \conv \{ u_{\tilde{H}}, u_{\hat{H}}\}$ and $\mathcal{E}_{\mathfrak{g}^{\ast}}(\mathfrak{a}_0) = \{X_{u}(\mathfrak{a}_0): u \in E^{\sharp}\} = \{ (1-t)X_{\tilde{u}}(\mathfrak{a}_0) + t X_{\hat{u}}(\mathfrak{a}_0), t\in [0,1] \}$. 
We claim that linear combinations of $X_{\tilde{u}}(\mathfrak{a}_0)$ and $X_{\hat{u}}(\mathfrak{a}_0)$ can't be tangent to $E$ whenever they are transversal to $E$. 
In fact, $X_{\tilde{u}}(\mathfrak{a}_0)$ and $X_{\hat{u}}(\mathfrak{a}_0)$ can be written as $a \tilde{n} + b\eta$ and $c \hat{n} + d \eta$, since they are parallel to $\tilde{H}$ and $\hat{H}$ respectively,
and $a,c\neq 0$ due to their transversality to $E$. 
Then $\left< X_{\tilde{u}}(\mathfrak{a}_0) \times X_{\hat{u}}(\mathfrak{a}_0), \eta \right> = \left< ac (\tilde{n} \times \hat{n}), \eta\right> \neq 0$, since $\tilde{n} \times \hat{n}$ is parallel to $E$, which proves that there isn't any linear combination of  $X_{\tilde{u}}(\mathfrak{a}_0)$ and $X_{\hat{u}}(\mathfrak{a}_0)$ which is tangent to $E$.
Therefore, we can't have $\mathfrak{a}([t_0,t_1]) \subset E$, which proves the lemma. 

\end{proof}

\section{Some locally non-minimizing extremals in \texorpdfstring{\(SE(2)\)}{SE(2)} and \texorpdfstring{\(\mathbb{H}\)}{H}.}
\label{in E(2) and H}

In this work we are interested in studying extremals that are one-parameter subgroups \(t \mapsto \exp(t u)\) of a Lie group \(G\), which arise in the following way:  \(\mathfrak{a}_0\) is a vertex of \(S_{F_*}\), \(H \subset S_F\) is the dual face of the vertex \(\mathfrak{a}_0\), and there exists a unique \( u_0 \in \operatorname{ri} H\) such that \(\mathfrak{a}_0([u_0, \cdot]) = 0\). 

In this case, \(\mathfrak{a}_0([u_0, \cdot]) = 0\) is equivalent to
\[
\mathfrak{a}([u, v]) = \mathfrak{a}(A(u \times v)) = (\mathfrak{a} \circ A)(u \times v) = ((\mathfrak{a} \circ A)^\sharp, u, v) = 0 \quad \text{for every } v \in \mathfrak{g},
\]
where \((\cdot, \cdot, \cdot)\) is the scalar triple product. {\em Therefore, \(\mathfrak{a}([u, \cdot]) = 0\) iff \((\mathfrak{a} \circ A)^\sharp\) and \(u\) are linearly dependent.} 
Since we are using Milnor bases on \(\mathfrak{g}\) and \(\mathfrak{g}^*\), then the coordinates of \((\mathfrak{a} \circ A)^\sharp\) are equal to the coordinates of \(\mathfrak{a} \circ A\).

For the sake of simplicity, we consider \(\{\mathfrak{a} = 1\} \coloneqq \{v \in \mathfrak{g} : \mathfrak{a}(v) = 1\}\), instead of \(\operatorname{ri} H\), in the forthcoming analysis.

\begin{lemma}
\label{lem:unique_u}
There exists a unique \(u \in \{\mathfrak{a} = 1\}\) such that \(\mathfrak{a}([u, \cdot]) = 0\) iff \(\mathfrak{a}((\mathfrak{a} \circ A)^\sharp) \neq 0\). In this case,
\[
u = \frac{1}{\mathfrak{a}((\mathfrak{a} \circ A)^\sharp)} (\mathfrak{a} \circ A)^\sharp.
\]
\end{lemma}

\begin{proof}
If \(\mathfrak{a} \circ A = 0\), then \(\mathfrak{a}([u, \cdot]) = 0\) for every \(u \in \{\mathfrak{a} = 1\}\).  
If \(\mathfrak{a} \circ A \neq 0\) and \(\mathfrak{a}((\mathfrak{a} \circ A)^\sharp) \neq 0\), then the unique \(u \in \{\mathfrak{a} = 1\}\) such that \(\mathfrak{a}([u, \cdot]) = 0\) is given by
\[
1 = \mathfrak{a}(u) = \mathfrak{a}(c (\mathfrak{a} \circ A)^\sharp) \quad \Rightarrow \quad c = \frac{1}{\mathfrak{a}((\mathfrak{a} \circ A)^\sharp)},
\]
and
\begin{align}
\label{align_unique_u}
u = \frac{1}{\mathfrak{a}((\mathfrak{a} \circ A)^\sharp)} (\mathfrak{a} \circ A)^\sharp.
\end{align}
Finally, if \(\mathfrak{a} \circ A \neq 0\) and \(\mathfrak{a}((\mathfrak{a} \circ A)^\sharp) = 0\), then \(\mathfrak{a}([u,\cdot])=0\) implies that \(u\) and \((\mathfrak{a}\circ A)^\sharp\) are linearly dependent, and there is no \(u\in\{\mathfrak{a}=1\}\) such that \(\mathfrak{a}([u,\cdot])=0\). Therefore, the condition \(\mathfrak{a}((\mathfrak{a} \circ A)^\sharp) \neq 0\) is equivalent to the existence of a unique \(u \in \{\mathfrak{a} = 1\}\) with \(\mathfrak{a}([u, \cdot]) = 0\), which is given by \eqref{align_unique_u}.
\end{proof}

The matrix \(A\) with respect to the Milnor bases corresponds to
\[
A = \begin{pmatrix}
0 & 0 & 0 \\
0 & \lambda & 0 \\
0 & 0 & 1
\end{pmatrix},
\]
with \(\lambda = 0\) for \(\mathfrak{h}\) and \(\lambda = 1\) for \(\mathfrak{e}(2)\) (see \cite{Milnor}).

From now on, fix
\(\mathfrak a_0=(\mathfrak a_1,\mathfrak a_2,\mathfrak a_3)\in\mathfrak g^*\)
and suppose that
\(
\mathfrak a_0\bigl((\mathfrak a_0\circ A)^\sharp\bigr)\neq0.
\)
If \(\mathfrak g=\mathfrak h\), then necessarily
\(\mathfrak a_3\neq0\), and Lemma~\ref{lem:unique_u} gives
\[
u_0=\left(0,0,\frac{1}{\mathfrak a_3}\right).
\]
If \(\mathfrak g=\mathfrak e(2)\), then
\((\mathfrak a_2)^2+(\mathfrak a_3)^2\neq0\) and
\begin{equation}\label{u E(2)}
u_0=
\left(
0,\,
\frac{\mathfrak a_2}{(\mathfrak a_2)^2+(\mathfrak a_3)^2},\,
\frac{\mathfrak a_3}{(\mathfrak a_2)^2+(\mathfrak a_3)^2}
\right).
\end{equation}

In the following examples, we construct an extremal in \(G\) by means of suitable concatenations of left cosets of one-parameter subgroups. Then we use this extremal to show that the extremal  $t \mapsto \exp (tu_0)$ corresponding to \(\mathfrak{a}\equiv\mathfrak{a}_0\) is not locally minimizing. Here, \(G\) denotes either of the three-dimensional matrix Lie groups \(SE(2)\) or \(\mathbb{H}\). Accordingly, whenever convenient, we represent elements of \(G\) and their Lie algebras either  by their matrix representations  or as ordered triples.

\subsection{The \texorpdfstring{\(SE(2)\)}{SE2} example.}

The matrix Lie group \(SE(2)\) is given by
\[
SE(2)
=
\left\{
\begin{bmatrix}
\cos x  & \sin x & z \\
-\sin x & \cos x & y \\
0       & 0      & 1
\end{bmatrix}
\;\middle|\;
x,y,z\in\mathbb{R}
\right\}.
\]

A basis of its Lie algebra \(\mathfrak{e}(2)\) is
\[
\begin{aligned}
e_1 &=
\begin{bmatrix}
0  & 1 & 0 \\
-1 & 0 & 0 \\
0  & 0 & 0
\end{bmatrix},
&
e_2 &=
\begin{bmatrix}
0 & 0 & 0 \\
0 & 0 & 1 \\
0 & 0 & 0
\end{bmatrix},
&
e_3 &=
\begin{bmatrix}
0 & 0 & 1 \\
0 & 0 & 0 \\
0 & 0 & 0
\end{bmatrix}.
\end{aligned}
\]

Consequently,
\[
\mathfrak{e}(2)
=
\left\{
a e_1+b e_2+c e_3
=
\begin{bmatrix}
0  & a & c \\
-a & 0 & b \\
0  & 0 & 0
\end{bmatrix}
\;\middle|\;
a,b,c\in\mathbb{R}
\right\}.
\]

 Its  Lie brackets are \([e_1,e_2]=e_3\), \([e_1,e_3]=-e_2\) and \([e_2,e_3]=0.\)

Let
\[
X
=
a e_1+b e_2+c e_3
=
\begin{bmatrix}
0  & a & c \\
-a & 0 & b \\
0  & 0 & 0
\end{bmatrix}
\in\mathfrak{e}(2).
\]

The one-parameter subgroup generated by \(X\) is the map \(t\longmapsto \exp(tX)\),
where
\[
\exp(tX)
=
\begin{cases}
\begin{bmatrix}
\cos(at)
&
\sin(at)
&
\dfrac{c\sin(at)+b\bigl(1-\cos(at)\bigr)}{a}
\\[1.2em]
-\sin(at)
&
\cos(at)
&
\dfrac{b\sin(at)+c\bigl(\cos(at)-1\bigr)}{a}
\\[1.2em]
0 & 0 & 1
\end{bmatrix},
& a\neq 0,
\\[4em]
\begin{bmatrix}
1 & 0 & tc \\
0 & 1 & tb \\
0 & 0 & 1
\end{bmatrix},
& a=0.
\end{cases}
\]

Consider an element
\[
g=
\begin{bmatrix}
\cos \bar{x}  & \sin \bar{x} & \bar{z} \\
-\sin \bar{x} & \cos \bar{x} & \bar{y} \\
0       & 0      & 1
\end{bmatrix}
\in SE(2).
\]

For
\[
h=
\begin{bmatrix}
\cos x  & \sin x & z \\
-\sin x & \cos x & y \\
0       & 0      & 1
\end{bmatrix}
\in SE(2),
\]
the left translation by \(g\) is defined by \(L_g(h)=gh\). Hence,

\[
\begin{aligned}
L_g(h)
&=
\begin{bmatrix}
\cos \bar{x}  & \sin \bar{x} & \bar{z} \\
-\sin \bar{x} & \cos \bar{x} & \bar{y} \\
0       & 0      & 1
\end{bmatrix}
\begin{bmatrix}
\cos x  & \sin x & z \\
-\sin x & \cos x & y \\
0       & 0      & 1
\end{bmatrix}
\\[1em]
&=
\begin{bmatrix}
\cos(\bar{x}+x)
&
\sin(\bar{x}+x)
&
y\sin \bar{x}+z\cos \bar{x}+\bar{z}
\\
-\sin(\bar{x}+x)
&
\cos(\bar{x}+x)
&
y\cos \bar{x}-z\sin \bar{x}+\bar{y}
\\
0 & 0 & 1
\end{bmatrix}.
\end{aligned}
\]

From now on, we present an example in \(SE(2)\) for which the
extremal  $t \mapsto \exp(tu_0)$ associated with   the constant Pontryagin extremal \((u(t),\mathfrak{a}(t)) \equiv (u_0, \mathfrak{a}_0)\) is not
locally minimizing. To this end, we construct a triangular face  of $S_F$ which is dual to \(\mathfrak{a}_0\).

We are interested in the case \( u_0 =(0,0,1)\). A general triangle $ABC$ of $S_F$, which contains $u_0$ in its relative interior, and such that the vertex $A$ lies in the $yz$-plane $\subset \mathfrak{g}$, is constructed in the following way: Consider the points \(A=(0,a^{2},a^{3})\) and \(B=(b^{1},b^{2},b^{3})\)
in \(\mathfrak{g}\), where \(a^{2} b^{1}\neq 0\). For \(\lambda_1,\lambda_2<0\), define \(C=
u_0+\lambda_1(A-u_0)+\lambda_2(B-u_0)\), that is,
\[
C= (c^1, c^2, c^3) = 
\bigl(
\lambda_2b^{1},\,
\lambda_1a^{2}+\lambda_2b^{2},\,
\lambda_1(a^{3}-1)+\lambda_2(b^{3}-1)+1
\bigr).
\]
We further assume that $a^2, b^1 > 0$, which suffices for our purposes.

Let \(\alpha=(x,y,z)\in\mathfrak{g}^*\). The vertices of the face
\(ABC\) determine the planes
\[
\Pi_V=\{\alpha:\langle V,\alpha\rangle=1\},\quad V\in\{A,B,C\},
\]
whose equations are
\[
a^{2}y+a^{3}z=1,\; b^{1}x+b^{2}y+b^{3}z=1,\;
\lambda_2b^{1}x+(\lambda_1a^{2}+\lambda_2b^{2})y+
[\lambda_1(a^{3}-1)+\lambda_2(b^{3}-1)+1]z=1.
\]

The planes \(\Pi_A\), \(\Pi_B\), and \(\Pi_C\)
contain the faces $A^{\ast}$, $B^{\ast}$ and $C^{\ast}$ in $S_{F_{\ast}}$, dual to the vertices $A,B,C \in S_{F}$ respectively. Their precise extent will not be relevant,
as our analysis is  local,  near their common point $\mathfrak{a}_0$, which is given by
\begin{equation}\label{mathfrak{a}}
\mathfrak{a}_{0}
=
\left(
\frac{a^{2}(1-b^{3})-b^{2}(1-a^{3})}{a^{2}b^{1}},
\frac{1-a^{3}}{a^{2}},
1
\right).
\end{equation}
From \ref{lem:unique_u} or \eqref{u E(2)}, \(\mathfrak{a}_0\) is maximized by
\[
u_0 = \left(
0,\,
\frac{a^{2}(1-a^{3})}{(1-a^{3})^2 + (a^{2})^2},\,
\frac{(a^{2})^2}{(1-a^{3})^2 + (a^{2})^2}
\right).
\]
 For our purposes, it is enough to choose the simplified case satisfying $a^{3}=1$, what implies that  \( u_0=(0,0,1)\). 
This choice yields 
\[
 \mathfrak{a}_0
=
\left(\frac{1-b^{3}}{b^{1}},0,1\right).
\]
The fact that $u_0$ is in the relative interior of a face of $S_F$, implies that $\mathfrak{a}_0$ is the unique point of $S_{F_{\ast}}$ that maximizes $u_0$, what implies that
\begin{equation}
\label{equation_Sfast_below_zleq1-2}
\{z=1\} \cap S_{F_\ast} = \{\mathfrak{a}_0\} \text{ and } S_{F_{\ast}} \subset \{z\leq 1\}.
\end{equation} 
It follows from \eqref{fildsbymatrix} that  $(u_0,\mathfrak{a}_0)$ is a constant Pontryagin extremal, what implies that its corresponding extremal is given by $\gamma_{\mathfrak{a}_0}(t) = \exp(tu_0)$.

We shall now
show that  $\gamma_{\mathfrak{a}_0}$  fails to be locally
minimizing.
To this end, observe that each face  $A^{\ast}, B^{\ast}$ and $C^{\ast}$ of 
\(S_{F_*}\subset\mathfrak{g}^*\),  containing \(\mathfrak{a}_0\),  has a maximizing control given by 
the vertices \(A\), \(B\), or \(C\) of the triangle \(ABC\) respectively. Therefore,
using \eqref{fildsbymatrix}, we can compute the extended geodesic field on
each face of \(S_{F_*}\).

\begin{figure}[ht]
    \centering
    \includegraphics[width=1\textwidth]{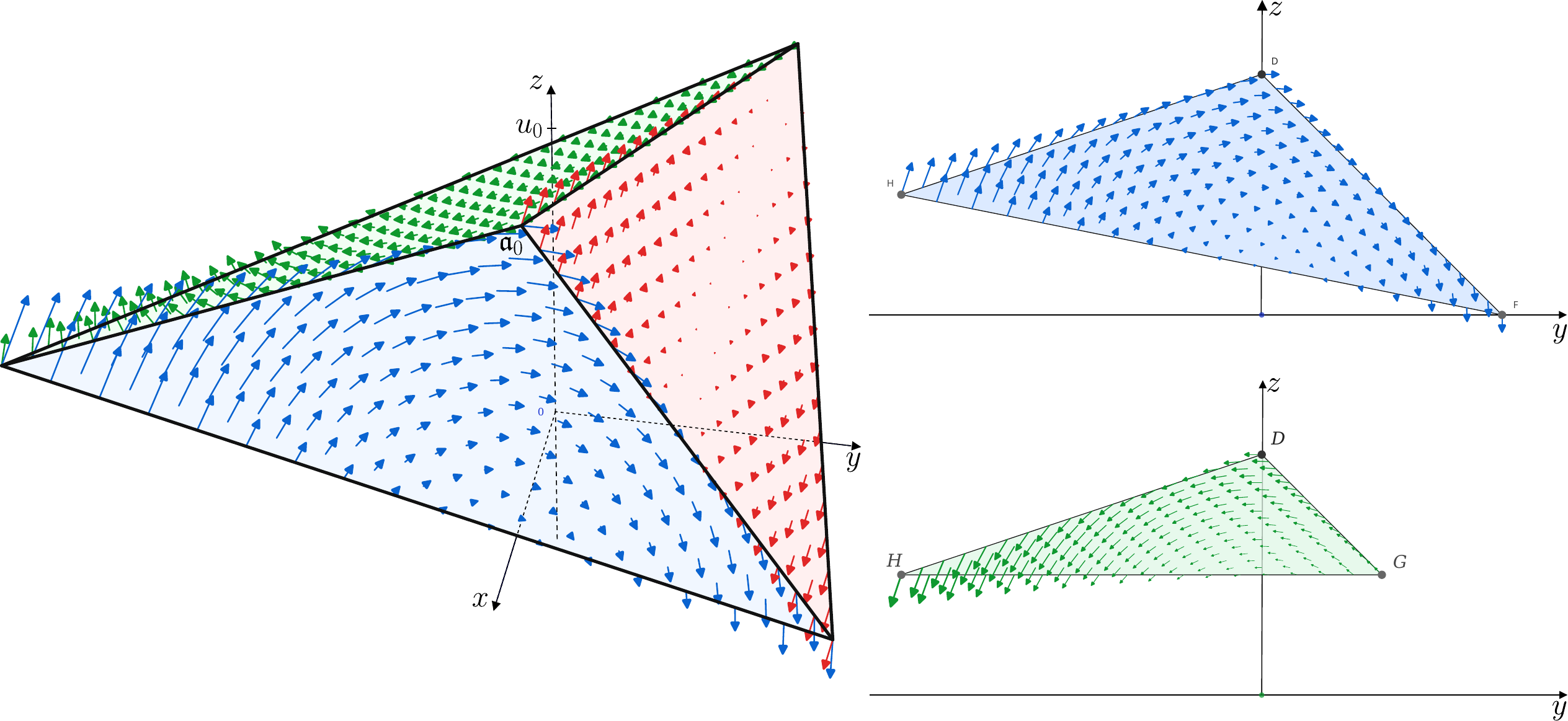}
    \caption{$\mathcal{E}_{\mathfrak{g}^*}$ on $S_{F_*}$}
    \label{ext geodesic on SF*}
\end{figure}

Note that, by Lemma \ref{across the edge} and Figure \ref{ext geodesic on SF*}, the integral curves
\(\mathfrak{a}(t)\) in a neighborhood of the vertex \(\mathfrak{a}_0\) cross
successive edges, thereby completing a full circuit around \(\mathfrak{a}_0\).
We now make this argument explicit.

For $SE(2)$, equation~\eqref{fildsbymatrix} takes the form
\[
X_u(\mathfrak{a})
=
\begin{bmatrix}
0 & u^3 & -u^2 \\
0 & 0 & u^1 \\
0 & -u^1 & 0
\end{bmatrix}
\begin{bmatrix}
\mathfrak{a}_1 \\
\mathfrak{a}_2 \\
\mathfrak{a}_3
\end{bmatrix}
=
\begin{bmatrix}
u^3\mathfrak{a}_2-u^2\mathfrak{a}_3 \\
u^1\mathfrak{a}_3 \\
-u^1\mathfrak{a}_2
\end{bmatrix}.
\]
Thus, $\mathcal{E}_{\mathfrak{g}^{\ast}}\vert_{\ri B^*}$ is given by
\[
X_B(\mathfrak{a})=
\begin{bmatrix}
b^3\mathfrak{a}_1-b^2\mathfrak{a}_3 \\
b^1\mathfrak{a}_3 \\
-b^1\mathfrak{a}_2
\end{bmatrix}.
\]
Let $\pi_{yz}: \mathfrak{g}^\ast \rightarrow Oyz$ be the orthogonal projection. Then $\pi_{yz}\vert_{B^{\ast}}$ is a diffeomorphism over its image and
\[
d\pi_{yz}(X_B(\mathfrak{a})) =  Y_{B}(\pi_{yz}(\mathfrak{a})) = b^1
\begin{bmatrix}
0 \\
\mathfrak{a}_3 \\
-\mathfrak{a}_2
\end{bmatrix},
\]
that is, $X_B$ and $Y_B$ are $\pi_{yz}$-related, what implies that $\alpha$ is an integral curve of $X_B$ iff $\pi_{yz}\circ \alpha$ is an integral curve of $Y_B$.
The integral curves of $Y_B$ are given by 
\begin{equation}
    \label{equation_curves_Bast}
\left(
0,
 \mathfrak{a}_2\cos(b^1t)+\mathfrak{a}_3\sin(b^1t),
 -\mathfrak{a}_2\sin(b^1t)+\mathfrak{a}_3\cos(b^1t)
\right).
\end{equation}
In the same way, $\mathcal{E}_{\mathfrak{g}^{\ast}}\vert_{\ri C^{\ast}}$ is given by 
\[
X_C(\mathfrak{a})=
\begin{bmatrix}
c^3\mathfrak{a}_1-c^2\mathfrak{a}_3 \\
c^1\mathfrak{a}_3 \\
-c^1\mathfrak{a}_2
\end{bmatrix}
\]
and its $\pi_{yz}$-related vector field on $\pi_{yz}(\ri C^*)$ is given by
\[
d\pi_{yz}(X_C(\mathfrak{a})) = Y_{C}(\pi_{yz}(\mathfrak{a})) = c^1
\begin{bmatrix}
0 \\
\mathfrak{a}_3 \\
-\mathfrak{a}_2
\end{bmatrix}.
\]
The integral curves of $Y_C$ are given by
\begin{equation}
    \label{equation_curves_Cast}
\left(
0,
 \mathfrak{a}_2\cos(c^1t)+\mathfrak{a}_3\sin(c^1t),
 -\mathfrak{a}_2\sin(c^1t)+\mathfrak{a}_3\cos(c^1t)
\right).
\end{equation}
Notice that $c^1 < 0$.

Finally, the face $A^{\ast}$ is parallel to the $x$-axis and $\mathcal{E}_{\mathfrak{g}^{\ast}}\vert_{\ri A^{\ast}}$ is given by
\[
X_A(\mathfrak{a})=
\begin{bmatrix}
\mathfrak{a}_2-a^2\mathfrak{a}_3 \\
0 \\
0
\end{bmatrix},
\]
and their integral curves are given by
\begin{equation}
    \label{equation_curva_Aast}
\big(\mathfrak{a}_1+t(\mathfrak{a}_2-a^2\mathfrak{a}_3),
     \mathfrak{a}_2,\mathfrak{a}_3\big).
\end{equation}

The vector fields $X_A$, $X_B$, $X_C$, $Y_A$, $Y_B$ and $Y_C$ will be continuously extended to the closure of their domains whenever we use Lemma \ref{across the edge}.

Firstly, we prove that the integral curves of $\mathcal{E}_{\mathfrak{g}^{\ast}}$ in a neighborhood of $\mathfrak{a}_0$ are closed.
All calculations below are straightforward (although sometimes lengthy), and we show only the results that are helpful for a better understanding. 
We will omit the $x$-coordinates when doing calculations in the $yz$-plane, whenever this does not cause confusion. 

The lines $\Pi_{A}\cap \Pi_{B}$, $\Pi_{A}\cap \Pi_{C}$ and $\Pi_{B}\cap \Pi_{C}$ contain the edges $A^{\ast} \cap B^{\ast}$, $A^{\ast} \cap C^{\ast}$ and $B^{\ast} \cap C^{\ast}$.
The projection $\rho_r := \pi_{yz}(\Pi_{A} \cap \Pi_{B})$ is given by $a^2 y + z = 1$ and so is $\pi_{yz}(\Pi_A \cap \Pi_C)$. 
The projection $\rho_l:=\pi_{yz}(\Pi_B \cap \Pi_C)$ is given by $a_l y + z = 1$, where $a_l = \frac{\lambda_1 a^2}{1 - \lambda_2} < 0$.
The projections of $A^{\ast} \cap B^{\ast}$ and $B^{\ast} \cap C^{\ast}$ on the $yz$-plane have $(0,1)$ as one of their vertices and the other vertices are placed below $(0,1)$ due to \eqref{equation_Sfast_below_zleq1-2}.
If $yz$-plane is oriented positively, then $A^{\ast} \cap B^{\ast}$ is placed on the right side of the $yz$-plane ($y\geq 0$, because $a^2 > 0$) and $B^{\ast} \cap C^{\ast}$ is placed on the left side ($y\leq 0$).

Consider the integral curves of $Y_B$ containing $(0,1-\varepsilon)$. 
They are  arc of circles oriented negatively (clockwise) (see \eqref{equation_curves_Bast}). 
Denote their trace by $C_{\varepsilon}:= S_{\Vert \cdot \Vert}(1-\varepsilon)$.
They intersect $\rho_r$ and $\rho_l$ at two points, whenever $\varepsilon \in (0,\varepsilon_0)$, where 
\[
\varepsilon_0 = \min \left\{ 1- \frac{1}{\sqrt{1+(a^2)^2}}, 1- \frac{1}{\sqrt{1+ a_l^2}}\right\}.
\]
Denote the intersection point of $\rho_{r} \cap C_{\varepsilon}$ which is nearer to $(0,1)$ by $q_{\varepsilon}$ and  the intersection point of $\rho_{l} \cap C_{\varepsilon}$ which is nearer to $(0,1)$ by $p_{\varepsilon}$.
Denote the point in $B^{\ast} \cap C^{\ast}$ that projects onto $p_{\varepsilon}$ by $P_{\varepsilon}$.
Likewise, denote the points in $A^{\ast} \cap B^{\ast}$ and $A^{\ast} \cap C^{\ast}$ that projects onto $q_{\varepsilon}$ by $Q^{AB}_{\varepsilon}$ and $Q^{AC}_\varepsilon$ respectively.
The points $Q^{AB}_{\varepsilon}$ and $Q^{AC}_\varepsilon$ lie in a line parallel to the $x$-axis, and the $x$-coordinate of $Q^{AB}_{\varepsilon}$ is greater than the $x$-coordinate of $Q^{AC}_{\varepsilon}$.

Since there exists an integral curve $\alpha_{B}(t)$ of $Y_B$ connecting $p_{\varepsilon}$ to $q_{\varepsilon}$, its corresponding integral curve $\beta_{B}(t)$ of $X_B$ on $B^{\ast}$ connects $P_{\varepsilon}$ to $Q^{AB}_{\varepsilon}$. 
The integral curve $\beta_B$ intersects $A^{\ast} \cap B^{\ast}$ at $Q^{AB}_{\varepsilon}$ transversally, because $\alpha_B$ intersects $\rho_r$ at $q_{\varepsilon}$ transversally.
From Lemma \ref{across the edge}, after reaching $Q^{AB}_\varepsilon$, the integral curve of $\mathcal{E}_{\mathfrak{g}^{\ast}}$ is immediately concatenated with the integral curve $\beta_A$ on $A^{\ast}$ (see \eqref{equation_curva_Aast}). 
Notice that it connects $Q^{AB}_{\varepsilon}$ to $Q^{AC}_{\varepsilon}$, because $\beta_A$ is parallel to the $x$-axis.
Of course, $\beta_A$ intersects $A^{\ast} \cap C^{\ast}$ at $Q^{AC}_{\varepsilon}$ transversally, because $A^{\ast} \cap C^{\ast}$ isn't parallel to the $x$-axis.
After reaching $Q^{AC}_{\varepsilon}$,  the integral curve of $\mathcal{E}_{\mathfrak{g}^{\ast}}$ is immediately concatenated with the integral curve $\beta_C$ on $C^{\ast}$ (again by Lemma \ref{across the edge}). 
Here, $q_{\varepsilon} = \pi_{yz}(Q^{AC}_{\varepsilon})$, the integral curve $\alpha_{C}$ of $Y_{C}$ beginning $q_{\varepsilon}$ has the same trace of $\alpha_{B}$, but oriented positively (see \eqref{equation_curves_Cast}).
Therefore it connects $q_{\varepsilon}$ to $p_{\varepsilon}$, and the corresponding curve $\beta_C$ on $C^{\ast}$ connects $Q^{AC}_{\varepsilon}$ to $P_{\varepsilon}$. 
Thus the integral curve of $\mathcal{E}_{\mathfrak{g}^{\ast}}$ is the concatenation $\beta_{\varepsilon} = \beta_B \cdot \beta_A \cdot \beta_C$, what is a closed curve.

Let us calculate the time $t_1(\varepsilon)$, $t_2(\varepsilon)$ and $t_3(\varepsilon)$, spent by $\beta_{\varepsilon}$ in each face $B^{\ast}$, $A^{\ast}$ and $C^{\ast}$ respectively. 
Recall that the vector fields $X_B$ and $X_C$ are $\pi_{yz}$-related to $Y_B$ and $Y_C$ respectively. 
Therefore, properties related to $X_B$ and $B^{\ast}$ can be studied through $Y_B$ and $\pi_{yz}(B^{\ast})$. 

The time $t_1(\varepsilon)$ is calculated dividing the angle elapsed from $p_{\varepsilon}$ to $q_{\varepsilon}$ by the angular velocity $b^1$ (see \eqref{equation_curves_Bast}). 
Straightforward, but lengthy calculation leads to
\begin{align}
t_1(\varepsilon)
& = 
\frac1{b^1}
\left[
\arctan a^2
-
\arccos\!\left(
\frac{1}{(1-\varepsilon)\sqrt{1+(a^2)^2}}
\right) 
\right] \nonumber \\
& +
\frac1{b^1}
\left[
\arctan\!\left(
\frac{-\lambda_1 a^2}{1-\lambda_2}
\right)
-
\arccos\!\left(
\frac{1}
{(1-\varepsilon)
\sqrt{
1+\left(\frac{\lambda_1 a^2}{1-\lambda_2}\right)^2
}}
\right)
\right]. \label{align_t1varepsilon}
\end{align}

The time $t_3(\varepsilon)$ is calculated in the same way. We get
\begin{equation}
    \label{equation_t3t1_fecha}
t_3(\varepsilon)=-\frac{t_1(\varepsilon)}{\lambda_2}.
\end{equation}

Finally, for $t_2(\varepsilon)$, we consider the Euclidean distance between $Q^{AB}_\varepsilon$ and $Q^{AC}_{\varepsilon}$ and divide it by the speed $\vert \mathfrak{a}^2 - a^2 \mathfrak{a}^3\vert$. We get
\[
t_2(\varepsilon)
=
\frac{a^2(1-\lambda_1-\lambda_2)}{-\lambda_2b^1}
\,
\frac{
a^2-\sqrt{\Delta_\varepsilon}
}{
(1+(a^2)^2)\sqrt{\Delta_\varepsilon}
}
\]
where
\[
\Delta_\varepsilon
=
(a^2)^2-(1+(a^2)^2)(2\varepsilon-\varepsilon^2).
\]
 
Henceforth, we consider the family of extremals in \(SE(2)\) connecting $e$ and the trace of $\gamma_{\mathfrak{a}_0}$, associated with the
successive passage of \(\mathfrak{a}(t)\) through the blue, red, and
green faces in Figure~\ref{ext geodesic on SF*}. It is obtained by
concatenation of left cosets of the exponentials of \(B\), \(A\), and \(C\), in this order. This concatenation will be used to show that $\gamma_{\mathfrak{a}_0}$ is not locally minimizing. 
Only after completing this step, we will return to the results established above about closed curves $\beta_\varepsilon$, to show that these concatenations are indeed Pontryagin extremals.

The one-parameter subgroups generated
by \(B\), \(A\), and \(C\) are given, respectively, by
\begin{equation}\label{gamma B}
\exp(tB)
=
\begin{pmatrix}
\cos(tb^{1})
&
\sin(tb^{1})
&
\dfrac{
b^{3}\sin(tb^{1})-b^{2}\cos(tb^{1})+b^{2}
}{b^{1}}
\\[2mm]
-\sin(tb^{1})
&
\cos(tb^{1})
&
\dfrac{
b^{2}\sin(tb^{1})+b^{3}\cos(tb^{1})-b^{3}
}{b^{1}}
\\[2mm]
0&0&1
\end{pmatrix},
\end{equation}
\begin{equation}\label{gamma A}
\exp(tA)
=
\begin{pmatrix}
1&0&t\\
0&1&ta^{2}\\
0&0&1
\end{pmatrix},
\end{equation}
and
\begin{equation}\label{gamma E}
\exp(tC)
=
\begin{pmatrix}
\cos(t\lambda_2b^{1})
&
\sin(t\lambda_2b^{1})
&
\dfrac{
\Delta\sin(t\lambda_2b^{1})
-\Lambda_2\cos(t\lambda_2b^{1})
+\Lambda_2
}{\lambda_2b^{1}}
\\[2mm]
-\sin(t\lambda_2b^{1})
&
\cos(t\lambda_2b^{1})
&
\dfrac{
\Lambda_2\sin(t\lambda_2b^{1})
+\Delta\cos(t\lambda_2b^{1})
-\Delta
}{\lambda_2b^{1}}
\\[2mm]
0&0&1
\end{pmatrix},
\end{equation}
where
\[
\Lambda_2=\lambda_1a^{2}+\lambda_2b^{2},
\qquad
\Delta=\lambda_2(b^{3}-1)+1.
\]

On the other hand, 
\begin{equation}\label{gamma a}
\gamma_{\mathfrak{a}_0}(t)=\exp(tu)
=
\begin{pmatrix}
1&0&t\\
0&1&0\\
0&0&1
\end{pmatrix}.
\end{equation}

By translations, the extremal obtained by concatenation can be written
as the piecewise smooth curve
\[
\gamma_{(t_1,t_2,t_3)}(s)=
\begin{cases}
\exp(sB), & 0\leq s\leq t_1,\\[1mm]
\exp(t_1B)\exp\bigl((s-t_1)A\bigr),
    & t_1\leq s\leq t_1+t_2,\\[1mm]
\exp(t_1B)\exp(t_2A)
\exp\bigl((s-t_1-t_2)C\bigr),
    & t_1+t_2\leq s\leq t_1+t_2+t_3.
\end{cases}
\]

To show that \(\gamma_{\mathfrak{a}_{0}}\) is not locally minimizing, it suffices
to prove that the competing extremal \(\gamma_{(t_1,t_2,t_3)}\) reaches the same point
as \(\gamma_{\mathfrak{a}_{0}}\) with strictly smaller length. To this end, we first
determine conditions on \(t_2\) and \(t_3\)  in terms of $t_1$ ensuring that \(\gamma\) reaches
the \(z\)-axis. We then analyze the evolution of its endpoint.

Set some variables to organize better some equations below:
\[
\theta=t_1b^{1},\qquad
\theta_2=t_3\lambda_2b^{1},\qquad
\Theta=\theta+\theta_2=b^{1}(t_1+\lambda_2t_3).
\]
A direct matrix multiplication yields
\[
\exp(t_1B)\exp(t_2A)\exp(t_3C)
=
\begin{pmatrix}
\cos\Theta & \sin\Theta & p(t_1)\\
-\sin\Theta & \cos\Theta & q(t_1)\\
0&0&1
\end{pmatrix},
\]
where
\begin{align}
p(t_1)
={}&
\frac{
b^{3}\sin\theta-b^{2}\cos\theta+b^{2}
}{b^{1}}
+t_2\bigl(\cos\theta+a^{2}\sin\theta\bigr)\nonumber
\\
&+
\frac{
\Delta\sin\Theta-\Lambda_2\cos\Theta
+\Lambda_2\cos\theta-\Delta\sin\theta
}{\lambda_2b^{1}},
\end{align}
\begin{align}
q(t_1)
={}&
\frac{
b^{2}\sin\theta+b^{3}\cos\theta-b^{3}
}{b^{1}}
+t_2\bigl(-\sin\theta+a^{2}\cos\theta\bigr)\nonumber
\\
&+
\frac{
\Lambda_2\sin\Theta+\Delta\cos\Theta
-\Lambda_2\sin\theta-\Delta\cos\theta
}{\lambda_2b^{1}},
\end{align}

To determine when the concatenated extremal meets the constant
extremal, write
\[
\gamma_{\mathfrak{a}_{0}}(T)
=
\begin{pmatrix}
1&0&T\\
0&1&0\\
0&0&1
\end{pmatrix}.
\]
Thus, the equality
\[
\exp(t_1B)\exp(t_2A)\exp(t_3C)=\gamma_{\mathfrak{a}_{0}}(T)
\]
is equivalent to
\[
\Theta\in 2\pi\mathbb{Z},
\qquad
q(t_1)=0,
\qquad
T=p(t_1).
\]
Since we are interested in $\gamma_{\mathfrak{a}_0}$ with values a neighborhood of the identity, we take
\(\Theta=0\). As
\[
\Theta=b^{1}(t_1+\lambda_2t_3),
\]
we obtain
\begin{equation}
    \label{equation_t3t1_eixoz}
t_3(t_1)=-\frac{t_1}{\lambda_2}.
\end{equation}
Substituting this expression into \(q(t_1)=0\), and setting
\(\theta=t_1b^{1}\),   yields
\[
t_2(t_1)
=
\frac{
a^{2}\lambda_1\sin\theta
-(1-\lambda_2)(1-\cos\theta)
}{
b^{1}\lambda_2
\bigl(a^{2}\cos\theta-\sin\theta\bigr)
}.
\]
Therefore,
\begin{equation}
    \label{equation_t2t1}
t_2(t_1)
=
\frac{
a^{2}\lambda_1\sin(t_1b^{1})
-(1-\lambda_2)\bigl(1-\cos(t_1b^{1})\bigr)
}{
b^{1}\lambda_2
\bigl(a^{2}\cos(t_1b^{1})-\sin(t_1b^{1})\bigr)
}
\end{equation}
provided that
\[
a^{2}\cos(t_1b^{1})-\sin(t_1b^{1})\neq 0,
\]
which holds for a sufficiently small $t_1$. 
This will be assumed henceforth.
In view of the dependence of \(t_2\) and \(t_3\) on \(t_1\), it is more
convenient to denote the parameter \(t_1\) simply by \(t\) from now on, and we
set  $\gamma_t:= \gamma_{(t,t_2(t),t_3(t))}$. 
Moreover, since we are interested in analyzing the evolution of the endpoint 
of \(\gamma_t\), we shall use \(\gamma\) 
to denote the curve describing its endpoint
\[
\gamma(t)=\exp(tB)\exp\bigl(t_2(t)A\bigr)\exp\bigl(t_3(t)C\bigr).
\] 
Thus
\[
\gamma(t)=\gamma_{\mathfrak{a}_{0}}\bigl(p(t)\bigr).
\]

Substituting
\[
t_3(t)=-\frac{t}{\lambda_2}
\]
and
\[
t_2(t)
=
\frac{
a^{2}\lambda_1\sin\theta
-(1-\lambda_2)(1-\cos\theta)
}{
b^{1}\lambda_2
\bigl(a^{2}\cos\theta-\sin\theta\bigr)
},
\qquad \theta=tb^{1},
\]
into
\[
p(t)
=
\frac{
a^{2}\lambda_1(\cos\theta-1)
+(\lambda_2-1)\sin\theta
}{
b^{1}\lambda_2
}
+t_2(t)\bigl(\cos\theta+a^{2}\sin\theta\bigr),
\]
we obtain

\begin{equation}
p(t)
=
\frac{
\bigl(1-\cos(tb^{1})\bigr)
\bigl((a^{2})^2\lambda_1+1-\lambda_2\bigr)
+a^{2}(\lambda_1+\lambda_2-1)\sin(tb^{1})
}{
b^{1}\lambda_2
\bigl(a^{2}\cos(tb^{1})-\sin(tb^{1})\bigr)
}.
\end{equation}
Since both \(\gamma_t\) and \(\gamma_{\mathfrak{a}_{0}}\) are parametrized by arc
length, their lengths up to the common endpoint
\[
\gamma(t)=\gamma_{\mathfrak{a}_{0}}\bigl(p(t)\bigr)
=
\begin{pmatrix}
1&0&p(t)\\
0&1&0\\
0&0&1
\end{pmatrix}
\]
are equal to the corresponding elapsed times. Therefore,
\[
\ell\left(\gamma_{\mathfrak{a}_{0}}|_{[0,p(t)]}\right)=p(t),
\qquad
\ell(\gamma_t)=t+t_2(t)+t_3(t).
\]
Consequently, to prove that \(\gamma_{\mathfrak{a}_{0}}\) is not locally
minimizing, it suffices to show that
\[
t+t_2(t)+t_3(t)<p(t)
\]
for all sufficiently small \(t>0\).

Let
\(
\ell_{\gamma_{\mathfrak{a}_{0}}}(t)=p(t)\) and
\(\ell_{\gamma}(t)=t+t_2(t)+t_3(t)
\)
denote the lengths of \(\gamma_{\mathfrak{a}_{0}}\) and \(\gamma_t\), respectively,
up to their common endpoint  $(0,0,p(t))$. 
Expanding the functions \(t_2\), \(t_3\),
and \(p\) in their Taylor series at \(t=0\), we obtain
\[
t_2(t)
=
\frac{\lambda_1}{\lambda_2}t
+
\frac{b^{1}(2\lambda_1+\lambda_2-1)}
     {2a^{2}\lambda_2}t^2
+O(t^3),
\qquad
t_3(t)=-\frac{t}{\lambda_2},
\]
and hence
\[
\ell_{\gamma}(t)
=
\frac{\lambda_1+\lambda_2-1}{\lambda_2}t
+
\frac{b^{1}(2\lambda_1+\lambda_2-1)}
     {2a^{2}\lambda_2}t^2
+O(t^3).
\]
On the other hand,
\[
\ell_{\gamma_{\mathfrak{a}_{0}}}(t)
=
\frac{\lambda_1+\lambda_2-1}{\lambda_2}t
+
\left[
\frac{b^{1}(2\lambda_1+\lambda_2-1)}
     {2a^{2}\lambda_2}
+
\frac{a^{2}b^{1}\lambda_1}{2\lambda_2}
\right]t^2
+O(t^3).
\]
Then
\[
\ell_{\gamma_{\mathfrak{a}_0}}(t)-\ell_{\gamma}(t) = \frac{a^{2} b^{1} \lambda_1}{2 \lambda_2}t^2 + O(t^3),
\]
what implies
\(
\ell_{\gamma_{\mathfrak{a}_{0}}}(t)>\ell_{\gamma}(t),
\)
or equivalently,
\(
p(t)>t+t_2(t)+t_3(t)
\)
for all sufficiently small \(t>0\). Hence, \(\gamma_t\) reaches the same
endpoint as \(\gamma_{\mathfrak{a}_{0}}\) with strictly smaller length, proving
that \(\gamma_{\mathfrak{a}_{0}}\) is not locally minimizing.

Now we prove that $\gamma_t$ is a Pontryagin extremal.
Firstly, for the closed integral curve $\beta_\varepsilon$ of $\mathcal{E}_{\mathfrak{g}^{\ast}}$, notice that $t_1(\varepsilon)$ is continuous and increasing with respect to $\varepsilon$ for a sufficiently small $\varepsilon$ because the angle $p_{\varepsilon} \hat{O} q_{\varepsilon}$ is continuous and increasing with respect to $\varepsilon$, and the angular velocity $b^1$ of $\beta_B$ does not depend on $\varepsilon$.
Likewise, the time interval $t_3(\varepsilon)$ is continuous and increasing with respect to $\varepsilon$.
The time interval $t_2(\varepsilon)$ is continuous and increasing with respect to $\varepsilon$ because the distance between $Q^{AB}_{\varepsilon}$ and $Q^{AC}_{\varepsilon}$ is continuous and increasing, and the modulus of the speed is continuous and decreasing with respect to $\varepsilon$ due to \eqref{equation_curva_Aast}.

Now consider $\eqref{equation_t2t1}$, and replace $t_1$ by $t_1(\varepsilon)$ given in \eqref{align_t1varepsilon} to compare $\gamma_t$ with a extremal corresponding to $\beta_{\varepsilon}$.
Then \eqref{equation_t2t1}, as a function of $\varepsilon$, is continuous and increasing for a sufficiently small $\varepsilon$, because it is the composition of two continuous and increasing functions. 
Denote it by $\bar t_2(\varepsilon)$.
Then the Taylor series of $\bar t_2(\varepsilon) - t_2(\varepsilon)$ centered at $\varepsilon = 0$ is given by
\begin{equation}
    \label{equation_t2bar_varepsilon_t2_varepsilon}
\bar t_2(\varepsilon) - t_2 (\varepsilon) = -\frac{(1-\lambda_2)(1-\lambda_1 - \lambda_2)^2}{\lambda_1^2 \lambda_2 b^1 (a^2)^3}\varepsilon^2 + O(\varepsilon^3).
\end{equation}
Notice that the coefficient of $\varepsilon^2$ is positive.

Then, for a sufficiently small positive value $T>0$, there exist $\bar{\varepsilon}_T$ and $\varepsilon_T$ such that $\bar t_2(\bar{\varepsilon}_T) = t_2(\varepsilon_T) = T$, where $\bar{\varepsilon}_T < \varepsilon_T$ due to \eqref{equation_t2bar_varepsilon_t2_varepsilon}.
Therefore, the time spent by $\gamma_t$ on $B^{\ast}$ when it spends $\bar t_2(\bar{\varepsilon}_T) = t_2(\varepsilon_T)=T$ on $A^{\ast}$ is $t_1(\bar{\varepsilon}_T) < t_1(\varepsilon_T)$, what also implies $t_3(\bar{\varepsilon}_T) < t_3(\varepsilon_T)$ due to \eqref{equation_t3t1_fecha} and \eqref{equation_t3t1_eixoz}. 
Hence, $\gamma_t$ is a proper restriction of the extremal corresponding to $\beta_{\varepsilon}$, what implies that $\gamma_t$ is an extremal for sufficiently small $t>0$.

\subsection{The Heisenberg  group.}

Consider the Heisenberg Lie group
\[
\mathbb{H}
=
\left\{
\begin{pmatrix}
1&x&z\\
0&1&y\\
0&0&1
\end{pmatrix}
:x,y,z\in\mathbb{R}
\right\},
\]
whose Lie algebra is
\[
\mathfrak{h}
=
\left\{
\begin{pmatrix}
0&a&c\\
0&0&b\\
0&0&0
\end{pmatrix}
:a,b,c\in\mathbb{R}
\right\}.
\]
We consider the Milnor basis
\[
e_1=
\begin{pmatrix}
0&1&0\\
0&0&0\\
0&0&0
\end{pmatrix},
\qquad
e_2=
\begin{pmatrix}
0&0&0\\
0&0&1\\
0&0&0
\end{pmatrix},
\qquad
e_3=
\begin{pmatrix}
0&0&1\\
0&0&0\\
0&0&0
\end{pmatrix},
\]
for which \([e_1,e_2]=e_3,\)
while all the remaining brackets vanish. In what follows, the elements of
\(\mathfrak{h}\) will be represented by their coordinates with respect
to the basis \(\{e_1,e_2,e_3\}\).

As in the \(SE(2)\) case, consider a polyhedral unit sphere in
\(\mathfrak h\) containing a triangular face \(ABC\), where
\(
A=(a^1,a^2,a^3),
\)
\(
B=(b^1,b^2,b^3),
\)
and assume
\(
\omega:=a^1b^2-a^2b^1>0.
\)
Since we are interested in the case where \(u_0=e_3=(0,0,1)\) is in the relative interior of $ABC$, we consider $\lambda_1, \lambda_2<0$ and \(C=u_0+\lambda_1(A-u_0)+\lambda_2(B-u_0).\)
Thus,
\(
C=(\Lambda^1,\Lambda^2,\Delta),
\)
where
\[
\Lambda^i=\lambda_1a^i+\lambda_2b^i,\qquad i=1,2,
\]
and
\[
\Delta
=
1+\lambda_1(a^3-1)+\lambda_2(b^3-1).
\]

Let \(\mathfrak a_0\in S_{F_*}\) be the vertex dual to \(ABC\). 
Since it is the unique point in $S_{F_{\ast}}$ maximizing $u_0$, then 
\begin{equation}
\label{equation_Sfast_below_zleq1}
\{z=1\} \cap S_{F_\ast} = \{\mathfrak{a}_0\} \text{ and } S_{F_{\ast}} \subset \{z\leq 1\}.
\end{equation} 

The faces
\(A^*,B^*,C^*\subset S_{F_*}\) containing \(\mathfrak a_0\) are maximized by
\(A,B,C\), respectively. We shall compare
\(
\gamma_{\mathfrak a_0}(t)=\exp(tu_0)
\) with a concatenation of left cosets of the one-parameter
subgroups generated by these controls.

For an element
\(
H=(h^1,h^2,h^3)\in\mathfrak{h},
\)
the associated one-parameter subgroup is
\begin{equation}
\exp(tH)
=
\begin{pmatrix}
1&th^1&th^3+\dfrac{t^2}{2}h^1h^2\\
0&1&th^2\\
0&0&1
\end{pmatrix}.
\label{eq:Heisenberg-exponential}
\end{equation}
Consequently,
\[
\exp(tA)
=
\begin{pmatrix}
1&ta^1&ta^3+\dfrac{t^2}{2}a^1a^2\\
0&1&ta^2\\
0&0&1
\end{pmatrix},
\]
\[
\exp(tB)
=
\begin{pmatrix}
1&tb^1&tb^3+\dfrac{t^2}{2}b^1b^2\\
0&1&tb^2\\
0&0&1
\end{pmatrix},
\]
and
\[
\exp(tC)
=
\begin{pmatrix}
1&t\Lambda^1&t\Delta+\dfrac{t^2}{2}\Lambda^1\Lambda^2\\
0&1&t\Lambda^2\\
0&0&1
\end{pmatrix}.
\]

Recall that the product in \(\mathbb{H}\) is given by
\[
\begin{pmatrix}
1&\bar{x}&\bar{z}\\
0&1&\bar{y}\\
0&0&1
\end{pmatrix}
\begin{pmatrix}
1&x&z\\
0&1&y\\
0&0&1
\end{pmatrix}
=
\begin{pmatrix}
1&\bar{x}+x&\bar{z}+z+\bar{x}y\\
0&1&\bar{y}+y\\
0&0&1
\end{pmatrix}.
\]

Now we show that integral curves of $\mathcal{E}_{\mathfrak{g}^{\ast}}$ in a neighborhood of $\mathfrak{a}_0$ are closed curves that successively traverses
the faces \(A^*\), \(B^*\), and \(C^*\). 
We describe their dynamics explicitly, as in the case of $\mathfrak{e}(2)$.
Firstly notice that there exist $\varepsilon \in (0,1)$ such that $z=1-\varepsilon$ intersects the relative interior of the edges $A^{\ast} \cap C^{\ast}$, $A^{\ast} \cap B^{\ast}$ and $B^{\ast} \cap C^{\ast}$ at $Q^{AC}_{\varepsilon}$, $Q^{AB}_{\varepsilon}$ and $Q^{BC}_{\varepsilon}$ respectively, due to \eqref{equation_Sfast_below_zleq1}.

\begin{figure}[ht]
    \centering
    \includegraphics[width=0.6\textwidth]{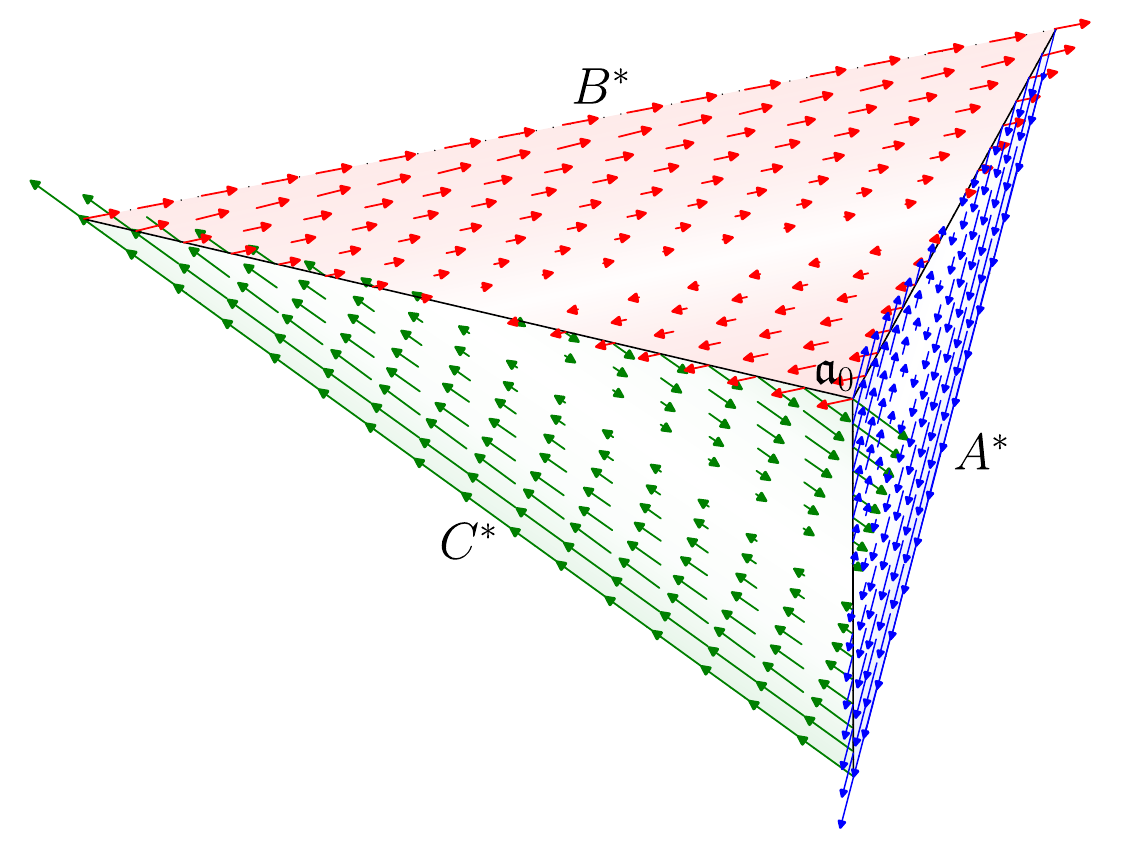}
    \caption{$\mathcal{E}_{\mathfrak{g}^*}$ on $S_{F_*}$}
\end{figure}

On $S_{F_{\ast}}$, the vector field along $A^*$,
as determined by \eqref{fildsbymatrix}, is given by
\begin{equation}
    \label{equation_XA}
\begin{aligned}
X_A(\mathfrak{a})
&=
\begin{pmatrix}
  0 & 0 & -a^2 \\
  0 & 0 &  a^1 \\
  0 & 0 & 0
\end{pmatrix}
\begin{pmatrix}
  \mathfrak{a}_1 \\
  \mathfrak{a}_2 \\
  \mathfrak{a}_3
\end{pmatrix}
=
\mathfrak{a}_3
\begin{pmatrix}
  -a^2 \\
   a^1 \\
   0
\end{pmatrix}.
\end{aligned}
\end{equation}
Consider $Q^{AC}_{\varepsilon}$ = ($A^{\ast} \cap C^{\ast}) \cap \{z=1-\varepsilon\}$ and set $(\mathfrak{a}_1^0, \mathfrak{a}_2^0, \mathfrak{a}_3^0) = Q^{AC}_{\varepsilon}$.
The integral curve of $X_A$ that begins at $Q^{AC}_{\varepsilon}$ (which is also an integral curve of $\mathcal{E}_{\mathfrak{g}^{\ast}}$) is
given by
\(
\beta_A(t) =\bigl(\mathfrak{a}^0_1-a^2\mathfrak{a}_3t,\,
      \mathfrak{a}_2^0+a^1\mathfrak{a}_3t,\,
      \mathfrak{a}_3^0\bigr) = \bigl(\mathfrak{a}^0_1-a^2(1-\varepsilon) t,\,
      \mathfrak{a}_2^0+a^1(1-\varepsilon) t,\,
      1-\varepsilon\bigr).
\)
Let us show that it reaches the face $B^*$ in positive time. 
To this end,
we determine the time $t$ for which

\[
  \left( (\mathfrak{a}_1^0,\mathfrak{a}_2^0,\mathfrak{a}_3^0)
  +(-a^2 (1-\varepsilon) t,\,a^1(1-\varepsilon) t,\,0)\right)(B)=1.
\]

Thus,
\[
t=
\frac{1-
 (\mathfrak{a}_1^0,\mathfrak{a}_2^0,\mathfrak{a}_3^0)(B)}
{(1-\varepsilon)(a^1b^2-a^2b^1)} > 0,
\]
because $B$ does not maximize $(\mathfrak{a}_1^0, \mathfrak{a}_2^0, \mathfrak{a}_3^0) \in A^{\ast} \cap C^{\ast}$ in $S_F$. 
The intersection point of $\beta_A$ with $A^{\ast} \cap B^{\ast}$ is $(\mathfrak{b}_1^0, \mathfrak{b}_2^0, \mathfrak{b}_3^0,):=Q^{AB}_{\varepsilon}$, because the $z$-coordinate of $\beta_A$ remains constant at $(1-\varepsilon)$.
Moreover, this intersection is transversal because $\beta_A$ has constant horizontal velocity and $A^{\ast} \cap B^{\ast}$ is not horizontal due to \eqref{equation_Sfast_below_zleq1}.
Therefore, due to Lemma \ref{across the edge}, the integral curve $\beta_A$ of $\mathcal{E}_{\mathfrak{g}^{\ast}}$ follows immediately to the face $B^{\ast}$.
But before going to the face $B^{\ast}$, we rewrite the time elapsed by $\beta_A$, which can be calculated straightforwardly, and it is given by
\begin{equation}
    \label{equation_tempoemAast_H}
    t_1(\varepsilon)=\frac{d_{\Vert \cdot \Vert}(Q^{AC}_{\varepsilon}, Q^{AB}_{\varepsilon})}{\Vert X_A(\mathfrak{a})\Vert} = \frac{\varepsilon (1-\lambda_1 - \lambda_2)}{-\lambda_2(1-\varepsilon)(a^1 b^2 - a^2 b^1)}.
\end{equation}

Once the trajectory reaches the face $B^*$, we similarly examine
whether the integral curve $\beta_B$ of $\mathcal{E}_{\mathfrak{g}^{\ast}}\vert_{B^\ast}$, which is the continuation of $\beta_A$, reaches the face $C$. Consider

\[
  \left( (\mathfrak{b}_1^0,\mathfrak{b}_2^0,\mathfrak{b}_3^0)
  +(-b^2(1-\varepsilon)t,\,b^1(1-\varepsilon)t,\,0)\right) (C)
=1,
\]

that is,
\[
t=
\frac{1-
  (\mathfrak{b}_1^0,\mathfrak{b}_2^0,\mathfrak{b}_3^0)(C)
}
{(1-\varepsilon)(b^1\Lambda^2-b^2\Lambda^1)}>0.
\]

Since $b^1\Lambda^2-b^2\Lambda^1
=\lambda_1(b^1a^2-b^2a^1)>0$ and $(\mathfrak{b}_1^0,\mathfrak{b}_2^0,\mathfrak{b}_3^0)(C)< 1$,
it follows that $t>0$. Hence, the trajectory  reaches $B^{\ast} \cap C^{\ast}$ transversally, and it follows immediately to the face $C^*$ due to Lemma \ref{across the edge}. The time interval that $\beta_B$ spends on $B^\ast$ is given by
\begin{equation}
    \label{equation_tempoemAast_H}
    t_2(\varepsilon)=\frac{d_{\Vert \cdot \Vert}(Q^{AB}_{\varepsilon}, Q^{BC}_{\varepsilon})}{\Vert X_B(\mathfrak{a})\Vert} = \frac{\varepsilon (1-\lambda_1 - \lambda_2)}{-\lambda_1(1-\varepsilon)(a^1 b^2 - a^2 b^1)}.
\end{equation}

Proceeding in the same way for the integral curve $\beta_C$ of $\mathcal{E}_{\mathfrak{g}^{\ast}}\vert_{C^{\ast}}$, we have that after spending
\begin{equation}
t(\varepsilon):=t_3(\varepsilon) = \frac{\varepsilon (1-\lambda_1 - \lambda_2)}{\lambda_1 \lambda_2 (1-\varepsilon) (a^1 b^2 -a^2 b^1)},
\end{equation}
on $C^{\ast}, \beta_C$ returns to the initial point $Q^{AC}_{\varepsilon}$. 
Therefore, the concatenation $\beta_{\varepsilon}:=\beta_A \cdot \beta_B \cdot \beta_C$ is a closed (periodic) integral curve of $\mathcal{E}_{\mathfrak{g}^{\ast}}$.
Notice that
\begin{equation}
    \label{equation_t1t2_em_funcao_de_t3}
    t_1(\varepsilon) = -\lambda_1 t_3(\varepsilon) \text{ and }t_2(\varepsilon) = -\lambda_2 t_3(\varepsilon).
\end{equation}

Now let us study concatenations of curves on $\mathbb{H}$, corresponding to curves $\mathfrak{a}(t)$ that
spends positive times \(t_1\), \(t_2\), and \(t_3\) on
\(A^*\), \(B^*\), and \(C^*\), respectively. 
Firstly we will not impose any restriction on $t_1$, $t_2$ and $t_3$, but we will see that the condition of reaching the orbit of $\gamma_{\mathfrak{a}_0}$, is equivalent to $t_1,t_2$ and $t_3$ satisfying \eqref{equation_t1t2_em_funcao_de_t3}

By left translations,
the corresponding curve, on \(\mathbb{H}\),  can be written as
\[
\gamma_{(t_1,t_2,t_3)}(s)
=
\begin{cases}
\exp(sA),
&0\leq s\leq t_1,\\[1mm]
\exp(t_1A)\exp\bigl((s-t_1)B\bigr),
&t_1\leq s\leq t_1+t_2,\\[1mm]
\exp(t_1A)\exp(t_2B)
\exp\bigl((s-t_1-t_2)C\bigr),
&t_1+t_2\leq s\leq t_1+t_2+t_3.
\end{cases}
\]
Analogously to the \(SE(2)\) example, we shall determine \(t_1\) and
\(t_2\) as functions of \(t_3\) and analyze the evolution of the endpoint
of \(\gamma\). In view of this similarity, we shall from the outset write simply
\[
\gamma(t)=\exp\bigl(t_1(t)A\bigr)\exp\bigl(t_2(t)B\bigr)\exp(tC).
\]
A direct matrix multiplication gives
\begin{equation}
\gamma(t)
=
\begin{pmatrix}
1&\xi(t)&\zeta(t)\\
0&1&\eta(t)\\
0&0&1
\end{pmatrix},
\label{eq:Heisenberg-concatenation}
\end{equation}
where
\[
\xi(t)
=
t_1a^1+t_2b^1+t\Lambda^1,
\qquad
\eta(t)
=
t_1a^2+t_2b^2+t\Lambda^2,
\]
and
\begin{align}
\zeta(t)
={}&
\frac{t_1^2}{2}a^1a^2+t_1a^3
+\frac{t_2^2}{2}b^1b^2+t_2b^3
+t_1t_2a^1b^2 \nonumber\\
&+
\frac{t^2}{2}\Lambda^1\Lambda^2
+t\Delta
+(t_1a^1+t_2b^1)t\Lambda^2.
\label{eq:Heisenberg-third-coordinate}
\end{align}

Recall that
\(
\mathfrak a_0(e_3)=1
\)
due to \eqref{equation_Sfast_below_zleq1} and that \(u_0=e_3\in\ri(ABC)\).
By Lemma~\ref{lem:unique_u}, the unique maximizing control associated with the
constant Pontryagin extremal \((u_0,\mathfrak a_0)\) is
\(u_0=e_3=(0,0,1)\). The corresponding one-parameter subgroup is therefore
\[
\gamma_{\mathfrak a_0}(\tau)
=
\exp(\tau e_3)
=
\begin{pmatrix}
1&0&\tau\\
0&1&0\\
0&0&1
\end{pmatrix}.
\]

We now determine when \(\gamma_{(t_1,t_2, t_3)}(t)\) reaches the trace of $\gamma_{\mathfrak{a}_0}$. By
\eqref{eq:Heisenberg-concatenation}, this occurs if and only if \(t_1a^1+t_2b^1
+t(\lambda_1a^1+\lambda_2b^1)=0\) and \(t_1a^2+t_2b^2
+t(\lambda_1a^2+\lambda_2b^2)=0\).
Since \(\omega=a^1b^2-a^2b^1>0\), the unique solution is
\begin{equation}
t_1(t)=-\lambda_1t,
\qquad
t_2(t)=-\lambda_2t,
\label{eq:Heisenberg-times}
\end{equation}
which are exactly the same conditions for $\mathfrak{a}(s)$ being a closed integral curve of $\mathcal{E}_{\mathfrak{g}^{\ast}}$ (see \eqref{equation_t1t2_em_funcao_de_t3}).
In particular, $\gamma_t := \gamma_{(t_1(t),t_2(t),t)}: [0,\tilde t] \rightarrow G$, where $\tilde t = t (1-\lambda_1 - \lambda_2)$, is an extremal for a sufficiently small $t>0$.

Let us see that $\gamma_{\mathfrak{a}_0}$ is not locally minimizing.
Substituting \eqref{eq:Heisenberg-times} into
\eqref{eq:Heisenberg-third-coordinate}, we obtain
\[
\gamma(t)
=
\exp\bigl(-\lambda_1tA\bigr)
\exp\bigl(-\lambda_2tB\bigr)
\exp(tC)
=
\begin{pmatrix}
1&0&\tau_{ABC}(t)\\
0&1&0\\
0&0&1
\end{pmatrix},
\]
where
\begin{equation}
\tau_{ABC}(t)
=
t(1-\lambda_1-\lambda_2)
+
\frac{\lambda_1\lambda_2}{2}
\omega t^2.
\label{eq:Heisenberg-ABE-endpoint}
\end{equation}
Since \(A\), \(B\), and \(C\) belong to \(S_F\), $\gamma_t$
is parametrized by arc length. Hence, its length up to the common endpoint $(0,0,\tau_{ABC}(t))$ is
\(
\ell(\gamma_t)
=
t(1-\lambda_1-\lambda_2).
\)
On the other hand, \(\gamma_{\mathfrak a_0}\) reaches the same endpoint at
time \(\tau_{ABC}(t) > (1-\lambda_1 - \lambda_2) t \).
Therefore, \(\gamma_t\) reaches the same
endpoint as \(\gamma_{\mathfrak{a}_{0}}\) with strictly smaller length, proving
that \(\gamma_{\mathfrak{a}_{0}}\) is not locally minimizing.

\begin{remark}
Let $G$ be a three-dimensional Lie group endowed with a left-invariant polyhedral Finsler structure $F$. 
Suppose that $v^1 \in S_F$ is in the relative interior of a face $H$ of $S_F$.
Consider a sequence of inner products $\tilde{g}_k$ in $\mathfrak{g}$ satisfying $v^1 \in S_{\tilde{g}_k}$ for every $k\in \mathbb{N}$, and such that the sequence of unit spheres $S_{\tilde{g}_k}$ converges to $S_F$ in a neighborhood of $v^1$ as $k \rightarrow \infty$. 
Extend the inner product $\tilde{g}_k$ to a left-invariant Riemannian structure $g_k$ on $G$.

Suppose that $t \mapsto \exp (tv^1)$ is an extremal of $(G,F)$, and let $\mathfrak{a}_0$ be the vertex of $S_{F_{\ast}}$, dual to $H$.
If $v_2 \in \mathfrak{g}$ is parallel to $H$ (that is, $v_2 \in \ker (\mathfrak{a}_0)$), then the sectional curvature of $\spann\{v_1,v_2\}$ in $(G,g_k)$ goes to infinity if $\ker \mathfrak{a}_0$ is not a subalgebra of $\mathfrak{g}$ (see Proposition 11 of \cite{Prudencio-Fukuoka-1}).
Here, it is important to observe that the sequence $g_k$ depends on $\mathfrak{a}_0$.
If we identify $v^1$ and $v^2$ with the corresponding left-invariant vector fields, and the lower bound of the sectional curvatures of $\spann \{v^1, v^2\}$ goes infinity, it makes sense to ask whether these extremals aren't locally minimizing (see Section \ref{section_introduction}). 
This section provides  two examples satisfying these properties.
\end{remark}

\bibliographystyle{amsplain}

\bibliography{references}{}

\end{document}